\input epsfx.tex
\overfullrule=0pt
 
\magnification=1400
\hsize=12.5cm
\vsize=16.5cm
\hoffset=-0.3cm   
\voffset=1.2cm    

\baselineskip 16 true pt

\font \caps=cmbx10 at 14 pt
\font \smcaps=cmbx10 at 12 pt
\font \pecaps=cmcsc10 at 9 pt

\newtoks \hautpagegauche  \hautpagegauche={\hfil}
\newtoks \hautpagedroite  \hautpagedroite={\hfil}
\newtoks \titregauche     \titregauche={\hfil}
\newtoks \titredroite     \titredroite={\hfil}
\newif \iftoppage        \toppagefalse   
\newif \ifbotpage        \botpagefalse    
\titregauche={\pecaps  Fran\c cois Dubois  }
\titredroite={\pecaps    Dual Raviart-Thomas mixed  finite elements }
\hautpagegauche = { \hfill \the \titregauche  \hfill  }
\hautpagedroite = { \hfill \the \titredroite  \hfill  }
\headline={ \vbox  { \line {  
\iftoppage    \ifodd  \pageno \the \hautpagedroite  \else \the
\hautpagegauche \fi \fi }     \bigskip  \bigskip  }}
\footline={ \vbox  {   \bigskip  \bigskip \line {  \ifbotpage  
\hfil {\oldstyle \folio} \hfil  \fi }}}

\def \smb {{\scriptstyle \bullet }}


\def\R{{\rm I}\! {\rm R}}

\def\ib#1{_{_{\scriptstyle{#1}}}}

\def\abs#1{\mid \! #1 \! \mid }

\def\mod#1{\setbox1=\hbox{\kern 3pt{#1}\kern 3pt}%
\dimen1=\ht1 \advance\dimen1 by 0.1pt \dimen2=\dp1 \advance\dimen2 by 0.1pt
\setbox1=\hbox{\vrule height\dimen1 depth\dimen2\box1\vrule}%
\advance\dimen1 by .1pt \ht1=\dimen1
\advance \dimen2 by .01pt \dp1=\dimen2 \box1 \relax}
\rm

\def\abs#1{\mid \! #1 \! \mid }

\def\mod#1{\setbox1=\hbox{\kern 3pt{#1}\kern 3pt}%
\dimen1=\ht1 \advance\dimen1 by 0.1pt \dimen2=\dp1 \advance\dimen2 by 0.1pt
\setbox1=\hbox{\vrule height\dimen1 depth\dimen2\box1\vrule}%
\advance\dimen1 by .1pt \ht1=\dimen1
\advance \dimen2 by .01pt \dp1=\dimen2 \box1 \relax}

\def\nor#1{\setbox1=\hbox{\kern 3pt{#1}\kern 3pt}%
\dimen1=\ht1 \advance\dimen1 by 0.1pt \dimen2=\dp1 \advance\dimen2 by 0.1pt
\setbox1=\hbox{\kern 1pt  \vrule \kern 2pt \vrule height\dimen1 depth\dimen2\box1
\vrule
\kern 2pt \vrule \kern 1pt  }%
\advance\dimen1 by .1pt \ht1=\dimen1
\advance \dimen2 by .01pt \dp1=\dimen2 \box1 \relax}

\def\sqr#1#2{{\vcenter{\vbox{\hrule height.#2pt \hbox{\vrule width .#2pt height#1pt 
\kern#1pt \vrule width.#2pt} \hrule height.#2pt}}}}
\def\square{\mathchoice\sqr64\sqr64\sqr{4.2}3\sqr33} 

\def\fl#1{\smash{\overrightarrow{#1}}}

\def\sth{\smash{^{^{\displaystyle \!\!\! \star}}}}


\hyphenation { li-mi-te vo-lu-mes vo-lu-me ma-the-ma-tik me-ca-nics me-thods 
ap-pro-xi-ma-tion 
}

$~$

\bigskip
\bigskip
\bigskip

\centerline{\caps  Dual Raviart-Thomas mixed  finite elements}

\bigskip
\bigskip
\bigskip

\centerline {\smcaps Fran\c{c}ois Dubois 
\footnote {$ ^{^{\displaystyle \smcaps \ast}}$}
{\rm  Conservatoire National des Arts et M\'etiers (Paris) et 
Universit\'e Paris-Sud (Orsay).}}

\bigskip
\centerline { April 08, 2002.
\footnote {$ ^{^{\displaystyle \square}}$}
{\rm Rapport  n$^{\rm o}$ 355  de l'Institut A\'ero-Technique de Saint Cyr, 
avril 2002.  Ce travail a \'et\'e pr\'esent\'e \`a  la troisi\`eme 
conf\'erence ``Finite Volumes for Complex Applications'', Porquerolles, juin 2002. 
Version abr\'eg\'ee publi\'ee avec le titre 
``Petrov-Galerkin Finite volumes'' dans le livre 
{\it  Finite Volumes for Complex Applications, Problems and Perspectives, volume 3}, 
Rapha\`ele Herbin et  Dietmar Kr\"oner Editeurs, Hermes Penton Science, pages  203-210, 2002. 
Edition du 7 novembre 2010.} }

\null\vskip 1cm  \noindent {\bf R\'esum\'e}

Pour un probl\`eme elliptique bidimensionnel, nous proposons de
formuler la m\'ethode des volumes   finis avec des \'el\'ements finis mixtes de
Petrov-Galerkin qui s'appuient sur la construction d'une base duale de
Raviart-Thomas.

\bigskip \bigskip \noindent {\bf Abstract}

For an elliptic problem with two space dimensions, we propose to formulate the finite
volume method with the help of Petrov-Galerkin mixed finite elementsthat are based on
the building of a dual Raviart-Thomas basis.


\bigskip \bigskip \noindent {\bf AMS (MOS) classification :} 65N30.
\null\vskip 1cm

\vfill\eject
\noindent {\bf Plan}
\smallskip

1)  Introduction

2)	 Stability analysis

3)	 Towards a first Petrov-Galerkin finite volume scheme

4)	 Perspectives

5)	 References.

\bigskip
\bigskip
\bigskip

\botpagetrue
\toppagetrue

\noindent {\smcaps 1)$\,\,\,$ Introduction.} 
\smallskip
\noindent $\bullet \quad$ Let $\, \Omega \,$ be a bidimensional bounded convex domain
in $\, \R^2 \,$ with a polygonal boundary $\,\, \partial \Omega. \,$ We consider   the
homogeneous Dirichlet problem for the Laplace operator in the domain $\, \Omega \,$~:

\smallskip \noindent  (1.1) $\qquad  \displaystyle
- \Delta u \,\,=\,\, f \,  \qquad{\rm in}\,\,    \Omega  \,$
\smallskip \noindent  (1.2) $\qquad  \displaystyle
\,\,\,\,  u \,\,\, \,\,\, \,\,=\,\, 0 \,  \qquad   {\rm on \,\, the \,\, boundary 
\,\, \partial \Omega  \,\, of} \,\, \Omega  . \,$

\smallskip \noindent
We suppose that the datum $\, f \,$ belongs to the space $\, L^2(\Omega) ,\,$ 

\smallskip \noindent  (1.3) $\qquad  \displaystyle
f \in L^2(\Omega)  \,, \, $

\smallskip \noindent
where this Hilbert  space  is  classically defined according to 

\setbox20=\hbox{$\displaystyle 
L^2(\Omega) \,\,= \,\, \Bigl\{ v \,:\, \Omega \longrightarrow \R \,,\,\, 
\int_{\Omega} \, \abs{v}^2 \, {\rm d}x \, < \infty \, \Bigr\} \,\, $ }
\setbox21=\hbox{$\displaystyle 
(u,\,v) \,\, \equiv \, \, \int_{\Omega} \,u \, v \,\, {\rm
d}x  \,\,,\qquad \forall \, u,\, v \,\, \in L^2(\Omega)  \,.  $ }
\setbox22=\hbox{$\displaystyle 
\parallel u \parallel\ib{0}^2 \, \,\,\equiv \,\,(u ,\,u ) \,\,,\qquad \qquad 
\forall \, u \,\, \in L^2(\Omega)  \,.  $ }
\setbox30= \vbox {\halign{#\cr \box20 \cr \box21 \cr \box22 \cr  }}
\setbox31= \hbox{ $\vcenter {\box30} $}
\setbox44=\hbox{\noindent (1.4) $\quad  \left\{ \box31 \right. \,$}  
\smallskip \noindent $ \box44 $

\smallskip \noindent
We introduce the momentum  $\,\, p \,\,$  defined by 

\smallskip \noindent  (1.5) $\qquad  \displaystyle
p \,\, = \,\, \nabla u \,.\,$

\smallskip \noindent
Taking  the divergence of both terms arising in equation (1.5), taking into
account the relation (1.1) and the hypothesis (1.3), we observe that 
the divergence of momentum $\,p \,$ belongs to the space $\, L^2(\Omega) .\,$ For this
reason, we introduce the vectorial Sobolev space $\, H({\rm div},\, \Omega) \,:\,$ 

\setbox20=\hbox{$\displaystyle  
H({\rm div},\, \Omega) \,\, = \,\, \bigl\{ \,q \in  L^2(\Omega) \times L^2(\Omega) \,
,\,\, {\rm div} \, q \in  L^2(\Omega) \, \bigr\} \,\, $ }
\setbox21=\hbox{$\displaystyle 
\parallel q \parallel  ^2 \ib{H({\rm div},\, \Omega)} \,\, = \,\, 
 \int_{\Omega} \, \bigl[ \, \abs{q}^2 \,+\, \abs{{\rm div}\,q}^2 \, \bigr] \, 
 {\rm d}x  \,\,,\qquad \forall \, q \,\, \in H({\rm div},\, \Omega) \,   $ }
\setbox30= \vbox {\halign{#\cr \box20 \cr \box21 \cr }}
\setbox31= \hbox{ $\vcenter {\box30} $}
\setbox44=\hbox{\noindent (1.6) $\quad  \left\{ \box31 \right. \,$}  
\smallskip \noindent $ \box44 $

\smallskip \noindent 
and we suppose in the following that the momentum $\, p \,$ satisfies  the condition 

\smallskip \noindent  (1.7) $\qquad  \displaystyle
p \, \in  H({\rm div},\, \Omega) \,.\,$ 

\bigskip \noindent $\bullet \quad$ 
The variational formulation of the  problem (1.1) (1.2) with the help of the pair $\,
\xi \,=\, (u ,\, p) \,$ is obtained by testing the definition (1.5) against a vector
valued function $ \, q \,$ and integrating by parts. With the help of the boundary
condition, it comes : 

\smallskip \noindent  (1.8) $\qquad  \displaystyle
(p,\,q) \,+\, (u,\, {\rm div} \, q) \,\,= \,\, 0 \,,\qquad \forall \, q \in 
 H({\rm div},\, \Omega) \,.  \,$ 

\smallskip \noindent  
Independently, the relations (1.1) and (1.5) are integrated on the domain $\, \Omega
\,$ after multiplying by a scalar valued function $\, v \, \in  L^2(\Omega) \,.\,$
We obtain : 

\smallskip \noindent  (1.9) $\qquad  \displaystyle
({\rm div} \, p ,\,v ) \,+\, (f,\, v) \,\,= \,\, 0 \,,\qquad \forall \, v \in 
L^2(\Omega) \,.  \,$ 

\smallskip \noindent
The ``mixed'' variational formulation is obtained by introducing the product space
$\, V \,$ defined as 

\setbox20=\hbox{$\displaystyle 
V \,\,= \,\, L^2(\Omega) \, \times \,  H({\rm div},\, \Omega) \,,\,$  }
\setbox21=\hbox{$\displaystyle 
\parallel (u,\,p)  \parallel\ib{V}^2 \,\,\,\,\equiv \,\,\, \parallel u
\parallel\ib{0}^2 \,\,+\,\, \parallel p \parallel\ib{0}^2 \,\,+\,\, \parallel {\rm
div}\,p  \parallel\ib{0}^2  \,\, \,,\,    $ }
\setbox30= \vbox {\halign{#\cr \box20 \cr \box21 \cr   }}
\setbox31= \hbox{ $\vcenter {\box30} $}
\setbox44=\hbox{\noindent (1.10) $\quad  \left\{ \box31 \right. \,$}  
\smallskip \noindent $ \box44 $

\smallskip \noindent 
the following bilinear form $\, \gamma(\smb,\, \smb) \,$ defined on $\, V \, \times \,
V \,: $ 

\smallskip \noindent  (1.11) $\qquad  \displaystyle
\gamma \bigl( (u,\,p),\, (v,\,q) \, \bigr) \,\,= \,\, (p,\,q) \,+\, (u,\, {\rm div}
\, q) \,+\,  ({\rm div} \, p ,\,v ) \,$ 

\smallskip \noindent 
and the linear form $\, \sigma(\smb) \,$ defined on $\, V   \, $ according to : 

\smallskip \noindent  (1.12) $\qquad  \displaystyle
< \sigma ,\, \zeta > \,\,=\,\, -(f,\,v) \,\,,\qquad \zeta = (v,\, q) \, \in V \,. \,$ 

\smallskip \noindent
Then the Dirichlet problem (1.1)(1.2) takes the form : 

\setbox20=\hbox{$\displaystyle 
\xi \, \in V  $ }
\setbox21=\hbox{$\displaystyle 
\gamma(\xi,\, \zeta) \,\,=\,\, < \sigma ,\, \zeta > \,\,,\qquad \forall \, \zeta \,
\in V \,.$ }
\setbox30= \vbox {\halign{#\cr \box20 \cr \box21 \cr }}
\setbox31= \hbox{ $\vcenter {\box30} $}
\setbox44=\hbox{\noindent (1.13) $\qquad  \left\{ \box31
\right. \,$}  
\smallskip \noindent $ \box44 $


\smallskip \noindent  
Due to classical inf-sup conditions introduced by Babu\v{s}ka [Ba71], the problem
(1.13) admits a unique solution $\, \xi \, \in V \,.\,$ 


\bigskip \noindent $\bullet \quad$ 
We introduce a mesh $\, {\cal T} \,$ that is   a bidimensional cellular complex
(see {\it e.g.} Godbillon [Go71])  composed in our case by triangular
elements $\, K \, $ $\, (K \, \in \,  {\cal E}\ib{\cal T} )  ,\,$  straight edges 
$\, a \, $ $\, (a \, \in \, {\cal A}\ib{\cal T}) \,$ and ponctual nodes $\, S \,$ 
$\, (S \, \in \, {\cal S}\ib{\cal T}) .\,$ We conside also classical finite
dimensional spaces $\, L\ib{\cal T}^2(\Omega) \,$ and $\,  H\ib{\cal T}({\rm div},\,
\Omega) \,$  that approximate  the spaces $\, L^2(\Omega) \,$ and $\,  H({\rm div},\,
\Omega)  \,$ respectively.  A scalar valued function $\, v \in  L_{\cal T}^2(\Omega)
\,$ is constant in each triangle $\, K \,$ of the mesh : 

\smallskip \noindent  (1.14) $\,\,\,  \displaystyle 
L\ib{\cal T}^2(\Omega) \,=\, \bigl\{ v : \Omega \longrightarrow \R ,\, \forall \, K 
\in {\cal E}\ib{\cal T} ,\, \exists \, v_K \, \in \R ,\, \forall \, x \in K, \, v(x)
\,=\, v_K \,\bigr\} \,.\,$

\smallskip \noindent 
A vector valued function  function $\, q \in  H\ib{\cal T}({\rm div},\, \Omega) \,$ is
a linear combination of Raviart-Thomas [RT77] basis functions $\, \varphi_a \,$ of
lower degree, defined in the forthcoming section. 

\bigskip \noindent $\bullet \quad$ 
Let $\, a \, \in {\cal A}\ib{\cal T} \,$ be an internal edge of the mesh, we denote by
$\, {\rm S} \,$ and $\, {\rm N} \,$ the two vertices that compose its boundary $\,
\partial a \,$  (see Figure 1) :

\smallskip \noindent  (1.15) $\quad  \displaystyle 
\partial a \,\,= \,\, \{ \,  {\rm S}  ,\, {\rm N} \,\} \,$ 

\smallskip \noindent
and by $\, K \,$ and $\,  L \,$ the two elements that compose its co-boundary $\,
\partial^c a \,$ 

\smallskip \noindent  (1.16) $\quad  \displaystyle 
\partial^c a \,\,= \,\, \{ \, K ,\, L \,\} \,$   

\smallskip \noindent
in such a way that the normal direction $\,n\ib{a}\,$ is oriented from $\, K \, $
towards $\, L \,$ and that the pair of vectors $\, (n\ib{a},\, \fl{ \rm SN}) \,$ is
direct, as shown on Figure 1. We denote by $\, {\rm W}\,$ (respectively by $\, {\rm E}
)\,$ the third vertex of the triangle $\, K \,$ (respectively of the triangle $\, L )
\,$ : 

\smallskip \noindent  (1.17) $\quad  \displaystyle 
K \,\,= \,\, ({\rm S},\, {\rm N} ,\, {\rm W} ) \,\,, \qquad 
L \,\,= \,\, ({\rm N},\, {\rm S} ,\, {\rm E} ) \,.\,$

\bigskip 
\centerline { \epsfysize=3cm    \epsfbox  {rt77.epsf} }
\smallskip  \smallskip

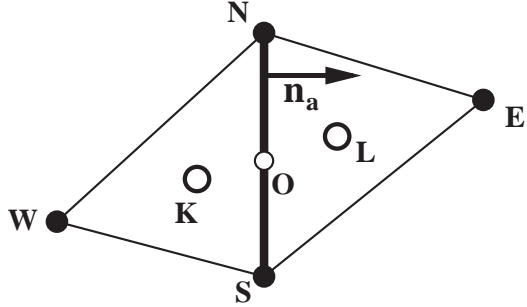
\centerline {    {\bf Figure 1.} Co-boundary $ \, (K,\,L) \,$ of the edge $\,
a=({\rm S},\,{\rm N}).\,$  } \smallskip

\smallskip \noindent
The vector valued   Raviart-Thomas [RT77] basis  function $\,\varphi_a \,$ is defined
by the relations 

\setbox20=\hbox{$\displaystyle 
\,\, {{1}\over{2 \abs{K}}} \, (x-{\rm W}) \,, \, \qquad x \, \in K  $ }
\setbox21=\hbox{$\displaystyle 
-{{1}\over{2 \abs{L}}} \, (x-{\rm E}) \,, \, \qquad x \, \in L  $ }
\setbox22=\hbox{$\displaystyle 
0 $ \qquad \qquad\qquad \qquad \quad   elsewhere.  }
\setbox30= \vbox {\halign{#\cr \box20 \cr \box21 \cr \box22 \cr}}
\setbox31= \hbox{ $\vcenter {\box30} $}
\setbox44=\hbox{\noindent (1.18) $\,\,\,  \varphi_a (x) \,\,= \,\, \left\{ \box31
\right. \,$}  
\smallskip \noindent $ \box44 $
 
\bigskip \noindent
When the edge $\, a \,$ is on the boundary $\, \partial \Omega \,$, we suppose that
the normal $\, n \,$ points towards the exterior of the domain, so the element $\, L
\,$ is absent. We have in all cases  the $\,  H({\rm div},\, \Omega) \,$ conformity : 

\smallskip \noindent  (1.19) $\quad  \displaystyle 
\varphi_a \, \in \,  H({\rm div},\, \Omega)  \,$

\smallskip \noindent
and the degrees of freedom are the fluxes of vector field $\, \varphi_a \,$ for all
the edges of the mesh (see [RT77]) : 

\smallskip \noindent  (1.20) $\quad  \displaystyle 
\int_b \, \varphi_a  \, \smb \, n\ib{a} \,\, {\rm d}\gamma \,\,= \,\, \delta_{a,\,b}
\,\,,\qquad \forall \, a ,\, b \, \in {\cal A}\ib{\cal T} \,.\, $ 

\smallskip \noindent
A vector valued function   $\, q \in  H\ib{\cal T}({\rm div},\, \Omega) \,$ is
a linear combination of the basis functions $\, \varphi_a \,: \,$ 

\smallskip \noindent  (1.21) $\quad  \displaystyle 
q \,=\, \sum_{a \, \in {\cal A}\ib{\cal T}} \, q_a \, \varphi_a \,  \,\,\,   \in
\, \,\,   H\ib{\cal T}({\rm div},\, \Omega) \,\,=\,\, < \varphi_b \,,\,\, b \, \in
{\cal A}\ib{\cal T} > \,.$ 

\bigskip \noindent $\bullet \quad$ 
The mixed finite element method consists in choosing as discrete linear space the
following product :

\smallskip \noindent  (1.22) $\quad  \displaystyle 
V\ib{\cal T} \,\,=\,\, L\ib{\cal T}^2(\Omega) \,\times \,  H\ib{\cal T}({\rm div},\,
\Omega) \,$ 

\smallskip \noindent
and proposes to replace the letter $\, V \,$ by $\, V\ib{\cal T}  \,$ inside the
variational formulation (1.13) :  

\setbox20=\hbox{$\displaystyle 
\xi\ib{\cal T}  \, \in V\ib{\cal T}   $ }
\setbox21=\hbox{$\displaystyle 
\gamma(\xi\ib{\cal T} ,\, \zeta) \,\,=\,\, < \sigma ,\, \zeta > \,\,,\qquad \forall \,
\zeta \, \in V\ib{\cal T}  \,$ }
\setbox30= \vbox {\halign{#\cr \box20 \cr \box21 \cr }}
\setbox31= \hbox{ $\vcenter {\box30} $}
\setbox44=\hbox{\noindent (1.23) $\qquad  \left\{ \box31
\right. \,$}  
\smallskip \noindent $ \box44 $

\smallskip \noindent
or in other terms 

\setbox20=\hbox{ $ \displaystyle   u\ib{\cal T} \, \in  L\ib{\cal T}^2(\Omega)
\,\,,\quad  p\ib{\cal T} \, \in   H\ib{\cal T}({\rm div},\, \Omega) \,$}
\setbox21=\hbox{ $ \displaystyle (p\ib{\cal T} \,,\, q ) \,+\,
(u\ib{\cal T} \,,\, {\rm div} \, q  ) \,\,=\,\, 0 \,,\qquad \forall \,
q \, \in   H\ib{\cal T}({\rm div},\, \Omega) \,$}
\setbox22=\hbox{ $ \displaystyle ({\rm div} \, p\ib{\cal T} \,,\, v )
\,+\, (f\,,\, v ) \,\, \,\,\,\,=\,\, 0 \,,\qquad \forall \, v  \,
\in  L\ib{\cal T}^2(\Omega) \,. \, $}
\setbox26= \vbox {\halign{# \cr \box20 \cr \box21 \cr  \box22 \cr }}
\setbox27= \hbox{ $\vcenter {\box26} $}
\setbox28=\hbox{\noindent (1.24) $\quad \left\{  \box27 \right. $}
\smallskip \noindent $ \box28 $

\smallskip \noindent  
The numerical analysis of the relations between the continuous problem (1.13) and the
discrete problem (1.23) as the mesh $\,{\cal T}\,$  is more and more refined is
classical [RT77]. The above method is popular in the context of petroleum and nuclear
industries   but suffers from the fact that the associated linear system is
quite difficult to solve from a practical point of view. The introduction of
supplementary Lagrange multipliers by Brezzi, Douglas and Marini [BDM85] allows a
simplification of these algebraic aspects,   their interpretation by Croisille in
the context of box schemes [Cr2k] gives a good mathematical foundation of a popular
numerical method and the possibility to reduce the size of the linear system 
has also been explored by Youn\`es, Mose, Ackerer and Chavent [YMAC97]. 

\bigskip \noindent $\bullet \quad$ 
From a theoretical and practical  point of view, the resolution of the system (1.24)
can be conducted as follows. We introduce the mass-matrix $\,\, M_{a,\,b} \, \,$
associated with the Raviart-Thomas vector valued functions : 

\smallskip \noindent  (1.25) $\quad  \displaystyle 
M_{a,\,b} \,\,=\,\, (\varphi_a ,\, \varphi_b) \,\,,\qquad a,\, b \, \in {\cal
A}\ib{\cal T} \,. \,$

\smallskip \noindent
Then the first equation of (1.24) determines the momentum 

\smallskip \noindent  (1.26) $\quad  \displaystyle 
p\ib{\cal T} \,\, \equiv \, \, \sum_{a \, \in   {\cal A}\ib{\cal T}} \,p\ib{{\cal
T}\!\!,a} \, \varphi_a \,\,$ 

\smallskip \noindent
as a function of the mean values $\,\, u\ib{{\cal T}\!\!,K} \, $ for $\, K \, \in 
{\cal E}\ib{\cal T} \,:\,$ 

\smallskip \noindent  (1.27) $\quad  \displaystyle 
p\ib{{\cal T}\!\!,a} \,\,=\,\,- \sum_{b \, \in  {\cal A}\ib{\cal T}}\,  \bigl( M^{-1}
\bigr)_{a,\,b}  \,\,  \sum_{ K \, \in  {\cal E}\ib{\cal T}} \, u\ib{{\cal T}\!\!,K}
\,\, \int_K \, {\rm div}\, \varphi_b \, {\rm d}x \,\,.   $ 

\smallskip \noindent 
The representation (1.27) suffers at our opinion form a major defect : due to the
fact that the matrix $\,\,  M^{-1} \,\,$ is full, the discrete gradient $
\,p\ib{\cal T} \,$ is a {\bf global} function of the mean values $\,  u\ib{{\cal
T}\!\!,K} \,$ and this property contradicts the mathematical  foundations of the
derivation operator to be  {\bf linear} and {\bf local}. An {\it a posteriori}
correction of this defect has been proposed by Baranger, Ma\^{\i}tre and Oudin [BMO96]
and with an appropriate  numerical integration of the mass matrix $\, M ,\,$ it is
possible to lump it and the discrete gradient in the direction $\, n\ib{a} \,$ of
the edge $\, a \,$ is represented by a formula of the type~: 

\smallskip \noindent  (1.28) $\quad  \displaystyle 
p\ib{{\cal T}\!\!,a} \,\,=\,\,{{u\ib{{\cal T}\!\!,L} - u\ib{{\cal T}\!\!,K}}
\over{h_a}} \,$

\smallskip \noindent
with the notations of Figure 1. The substitution of the relation (1.28) inside the
second equation of the formulation (1.24) conducts to a variant of the so-called
finite volume method. In an analogous manner, the family of finite volume schemes
proposed by Herbin [He95] supposes {\it a priori} that the discrete gradient in the
normal direction admits a representation of the form (1.28). Nevertheless, the
engineer  intuition is not correctly satisfied by a scheme such that (1.28). The finite
difference $\,\, \smash{{{u_{{\cal T}\!\!,L} - u_{{\cal T}\!\!,K}}\over{h_a}}} \,\,$
must be {\it a priori} to be a good (strong~?) approximation of the gradient $\,\,
p\ib{\cal T} \,=\, \nabla u\ib{\cal T} \,\,$ in the direction $\fl{KL} $ whereas the
coefficient $\,\, p\ib{{\cal T}\!\!,a} \,\,$ is an approximation of $\,\,  \int_{\rm
\displaystyle a} \nabla  u\ib{\cal T} \smb n\,\, {\rm d}\tau \,\,$ in the {\bf normal}
direction (see again the Figure 1). When the mesh $\, {\cal T} \,$ is composed by
general triangles, this approximation is not completely satisfactory at our opinion
and contains  a real limitation of these variants of the finite volume method. 

\bigskip \noindent $\bullet \quad$ 
We recall here that  the finite volume method for the approximation of the diffusion
operators  has been first proposed from empirical considerations. Following {\it e.g.}
Noh [No64] and Patankar [Pa80], the idea is to represent the normal interface gradient 
 $\,\,  \int_{\rm \displaystyle a}   \nabla  u\ib{\cal T} \smb n\,\, {\rm d}  \tau
\,\,$ as a function of {\bf neighbouring} values. Given an edge $\, a,\,$ a  vicinity
$\, {\cal V}(a) \,$ is {\bf first} determined in order to represent the normal gradient
$\,\,p\ib{{\cal T}\!\!,a}  =  \int_{\rm \displaystyle a}   \nabla  u\ib{\cal T} \smb
n\,\, {\rm d} \tau \,\,$  with a ``derivation formula''  of the type 

\smallskip \noindent  (1.29) $\quad  \displaystyle 
 \int_{\rm \displaystyle a}   \nabla  u\ib{\cal T} \smb n\,\, {\rm d} \tau \,\,=\,\,
 \sum_{K  \in {\cal V}(a)} \, g\ib{a,K} \,\, u\ib{{\cal T}\!\!,K}  \,. \,$ 
 
\smallskip \noindent  
Then the conservation equation 

\smallskip \noindent  (1.30) $\quad  \displaystyle 
{\rm div} \, p \,\,+ \,\, f \,\, = \,\, 0 \,$ 

\smallskip \noindent
is integrated inside each cell $\, K  \in {\cal E}\ib{\cal T} \,$ is order to
determine an equation for the mean values $\,  u\ib{{\cal T}\!\!,K} \,\,$ for all $\,
K  \in {\cal E}\ib{\cal T} .\,$ The difficulties of such approches have been presented
by Kershaw [Ke81] and a variant of such  scheme has been   first analysed by
Coudi\`ere, Vila and Villedieu [CVV99]. The key remark that we have done with
F.~Arnoux (see [Du89]), also observed by Faille, Gallou\"et and Herbin [FGH91]  is
that the representation (1.29) must be {\bf exact} for linear functions $\, u\ib{\cal
T} \,.\,$ We took this remark as a starting point for our tridimensional 
finite volume scheme proposed in [Du92]. It is also an essential hypothesis for the
result proposed by Coudi\`ere, Vila and Villedieu.

\bigskip \noindent $\bullet \quad$ 
In this contribution, we propose to  discretize the variational problem (1.13) with
the Petrov-Galerkin mixed finite element method, first introduced by Thomas and
Trujillo [TT99]. In the way  we have proposed  in [Du2k], the idea is to
construct a  discrete  functional  space 
$\,\,  H\ib{\cal T}\sth({\rm div},\, \Omega) \, \,$   generated by vectorial functions 
$\, \varphi^{\star}_{a} ,\,  a   \in {\cal A}\ib{\cal T} ,\,$ 
that are conforming in the space $\, H({\rm div},\, \Omega) \,$  

\smallskip \noindent  (1.31) $\quad  \displaystyle 
\varphi^{\star}_{a} \,\, \in \, H({\rm div},\, \Omega) \,$

\smallskip \noindent 
and to represent exactly  the {\bf dual basis} of the family $\, \{ \varphi_{b},\, b 
\in {\cal A}\ib{\cal T} \}\,$  with the $\, L^2 \,$ scalar product~: 

\smallskip \noindent  (1.32) $\quad  \displaystyle 
( \, \varphi_{a} \,,\,\varphi_{b}^{\star} \, ) \,\,= \,\, \delta_{a,\,b}
\,\,,\qquad \forall \, a ,\, b \, \in {\cal A}\ib{\cal T} \,.\, $ 

\smallskip \noindent  (1.33) $\quad  \displaystyle 
H\ib{\cal T}\sth({\rm div},\, \Omega) \,\,=\,\, <
\varphi_{b}^{\star}  \,,\,\, b \, \in {\cal A}\ib{\cal T} > \,\,\, \subset \,\,\, 
 H({\rm div},\, \Omega) \,.$ 

\smallskip \noindent
Then the mixed Petrov-Galerkin mixed finite element method consists just in
replacing  the space $\,\, H\ib{\cal T} ({\rm div},\, \Omega) \,\,$ by the dual space 
$\,\, H\ib{\cal T}\sth({\rm div},\, \Omega) \,\,$  for {\bf test functions} in the
first equation of discrete formulation (1.24). We obtain by doing this the so-called
{\bf  Petrov-Galerkin finite volume} scheme : 

\setbox20=\hbox{ $ \displaystyle   u\ib{\cal T} \, \in  L\ib{\cal T}^2(\Omega)
\,\,,\quad  p\ib{\cal T} \, \in   H\ib{\cal T}({\rm div},\, \Omega) \,$}
\setbox21=\hbox{ $ \displaystyle (p\ib{\cal T} \,,\, q ) \,+\,
(u\ib{\cal T} \,,\, {\rm div} \, q  ) \,\,=\,\, 0 \,,\qquad \forall \,
q \, \in    H\ib{\cal T}\sth({\rm div},\, \Omega)  \,$}
\setbox22=\hbox{ $ \displaystyle ({\rm div} \, p\ib{\cal T} \,,\, v )
\,+\, (f\,,\, v ) \,\, \,\,\,\,=\,\, 0 \,,\qquad \forall \, v  \,
\in  L\ib{\cal T}^2(\Omega) \,. \, $}
\setbox26= \vbox {\halign{# \cr \box20 \cr \box21 \cr  \box22 \cr }}
\setbox27= \hbox{ $\vcenter {\box26} $}
\setbox28=\hbox{\noindent (1.34) $\quad \left\{  \box27 \right. $}
\smallskip \noindent $ \box28 $

\smallskip \noindent
We introduce a compact form of the previous mixed Petrov-Galerkin formulation with the
help of the product space $\, V\ib{\cal T}^{\displaystyle \star} \,$ defined by 

\smallskip \noindent  (1.35) $\quad  \displaystyle 
V\ib{\cal T}^{\displaystyle \star} \,\,=\,\, L\ib{\cal T}^2(\Omega) \,\times \, 
H\ib{\cal T}\sth ({\rm div},\, \Omega) \,. \,$ 

\smallskip \noindent 
Then the formulation (1.34) admits the form : 

\setbox20=\hbox{$\displaystyle 
\xi\ib{\cal T}  \, \in V\ib{\cal T}   $ }
\setbox21=\hbox{$\displaystyle 
\gamma(\xi\ib{\cal T} ,\, \zeta) \,\,=\,\, < \sigma ,\, \zeta > \,\,,\qquad \forall \,
\zeta \, \in V\ib{\cal T}^{\displaystyle \star}  \,. \,$ }
\setbox30= \vbox {\halign{#\cr \box20 \cr \box21 \cr }}
\setbox31= \hbox{ $\vcenter {\box30} $}
\setbox44=\hbox{\noindent (1.36) $\qquad  \left\{ \box31
\right. \,$}  
\smallskip \noindent $ \box44 $

\smallskip \noindent
By doing this choice, it is easy to check that the scheme (1.34) is in fact a finite
volume scheme for the Laplace operator. The key point is to construct the so-called
{\bf  dual Raviart-Thomas basis functions } $\, \varphi_{a}^{\star}\,$ in order to
guaranty  Babu\v{s}ka's [Ba71] inf-sup stability property.

\bigskip \noindent $\bullet \quad$ 
The plan of the article is the following : we  derive in the second part sufficient
conditions in order to guaranty the final stability of the finite element scheme. Then
we propose a particular family of dual Raviart-Thomas functions and propose by doing
this a two-parameter family of finite volumes schemes. 


\bigskip
\bigskip
\noindent {\smcaps 2)$\,\,\,$ Stability analysis} 
\smallskip
\noindent $\bullet \quad$ 
We suppose in the following that the mesh $\, \, {\cal T} \,\,$ is a bidimensional 
cellular complex composed by triangles as
proposed in the first section.  Following the work of Ciarlet and Raviart [CR72], for
any element $\,\, K   \in \,  {\cal E}\ib{\cal T} \, \,$ we denote by $\,h\ib{\!K}\,$
the diameter of the triangle $\, K \,$ and by $\, \rho\ib{\!K}\,$ the diameter of the
inscripted ball inside $\, K .\,$ We suppose that the mesh $\, \, {\cal T} \,\,$
belongs   to a family $\, {\cal U}\ib{\theta} \,$ of meshes that satisfies the
following definition. 

\bigskip \noindent
{\bf Definition 1. \quad Family of regular meshes}

\noindent
Let $\, \theta \,$ be a strictly positive parameter. The family $\,
{\cal U}\ib{\theta} \,$ of meshes is defined by the condition 

\noindent  (2.1) $\quad  \displaystyle 
{\cal T} \,\in \, {\cal U}\ib{\theta} \,\,\, \Longleftrightarrow \,\,\, 
\forall \, K \, \in {\cal E}\ib{\cal T} \,,\,\, {{h\ib{\!K}}\over{\rho\ib{\!K}}} \,
\leq \, \theta \,.$

\smallskip \noindent
We suppose also that the dual space $\,\,  H\ib{\cal T}\sth({\rm div},\, \Omega) \,\,$
constructed by the conditions (1.31), (1.32), (1.33) satisfies the following
hypothesis. 

\bigskip \noindent
{\bf Hypothesis 1. \quad Interpolation operator $\,\, \,  H\ib{\cal T}({\rm div},\,
\Omega) \, \longrightarrow \, H\ib{\cal T}\sth ({\rm div},\,
\Omega) \,. \,$ }

\noindent
We suppose that the mesh $\, {\cal T} \,$ belongs to the family $\, {\cal
U}\ib{\theta} \,$ of Definition 1 and that the dual basis $\,\, \varphi^{\displaystyle 
\star}_a  \,\,$  is constructed in such a way that there exists a linear mapping 
 $\,\, H\ib{\cal T}({\rm div},\, \Omega) \, \ni
\, q \,\, \longmapsto \,\, \Pi \, q \,\in \, H\ib{\cal T}\sth ({\rm div},\, \Omega)
\,\,$ and strictly positive constants $\, A,\, B,\,
D  ,\, E \,$ that only depends on the parameter $\, \theta \,$ such that we have the
following estimations : 

\smallskip \noindent  (2.2) $\quad  \displaystyle 
A \, \parallel q \parallel\ib{0}^2 \,\,\,\, \leq \,\, ( \, q  \, ,\, \Pi \, q \, )
\,\,, \qquad \forall \, q \, \in  H\ib{\cal T}({\rm div},\, \Omega) \,$ 
\smallskip \noindent  (2.3) $\quad  \displaystyle 
\parallel \Pi \, q \parallel\ib{0} \,\,\,\,\, \leq \,\, B \, \parallel  q
\parallel\ib{0} \,\,, \qquad \forall \, q \, \in  H\ib{\cal T}({\rm div},\, \Omega)
\,$ 
\smallskip \noindent  (2.4) $\quad  \displaystyle 
\parallel {\rm div} \,\Pi \, q \parallel\ib{0}  \,\,\,\, \leq \,\, D \,
\parallel {\rm div} \, q \parallel\ib{0}  \,\,, \qquad \qquad 
\forall \, q \, \in  H\ib{\cal T}({\rm div},\, \Omega) \,$ 
\smallskip \noindent  (2.5) $\quad  \displaystyle 
( \, {\rm div}\, q \,,\,  {\rm div}\, \Pi \,q \,) \,\,\, \geq \,\, E \, 
\parallel {\rm div} \, q \parallel\ib{0}^2 \,\,, \qquad \forall \, q \, \in  H\ib{\cal
T}({\rm div},\, \Omega) \,. \,$ 

\bigskip 
\noindent
{\bf Proposition 1. \quad Divergence lifting of scalar fields}

\noindent
Let $\, \theta \,$ be a strictly positive parameter. We suppose that the dual
Raviart-Thomas basis satisfies the Hypothesis 1. Then there exists some strictly
positive constant $\, F \,$ that only depends on the parameter $\, \theta \,$ such that
for any mesh   $\, {\cal T} \,$  that belongs to the family $\, {\cal U}\ib{\theta} ,
\,$ and for any scalar field $\, u \,$ constant in each element $\, K  \in {\cal
E}\ib{\cal T} \,$ $(u \, \in \, L\ib{\cal T}^2(\Omega) ) ,\, $
there exists some vector field $\,\, q \, \in  \, H\ib{\cal T}\sth ({\rm div},\,
\Omega) \,\,$ such that 

\smallskip \noindent  (2.6) $\quad  \displaystyle 
\parallel q \parallel  \ib{H({\rm div},\, \Omega)} \,\,\, \leq \,\,\, F \, \parallel
u \parallel \ib{0} \,$ 

\smallskip \noindent  (2.7) $\quad  \displaystyle 
(\, u \,,\, {\rm div} \, q \, ) \quad \,  \geq \,\, \,\, \parallel u \parallel \ib{0}
^2\, \,. \, $ 

\smallskip \noindent
{\bf Proof of Proposition 1}. 

\smallskip \noindent $\bullet \quad$ 
Let $\,\, u  \in \, L\ib{\cal T}^2(\Omega)  \, $ be a discrete scalar function
supposed to be   constant in each triangle $\, K \,$ of the mesh $\,  {\cal T} . \,$
Let $\, \psi \, \in H^1_0(\Omega) \,$ be the variational solution of the Poisson
problem 

\smallskip \noindent  (2.8) $\quad  \displaystyle 
\Delta \, \psi \,\, = \,\, u \quad  {\rm in} \,\, \Omega \,\,,\qquad \qquad
\psi \,=\, 0 \quad   {\rm on} \,\, \partial \Omega \,.\,$ 

\smallskip \noindent  
Since $\, \Omega \,$ is convex, the solution $\, \psi \,$ of the problem (2.8) belongs
to the space $\, H^2(\Omega) \,$ and there exists some constant $\, G > 0 \,$ that only
depends on $\, \Omega \,$ such that 

\smallskip \noindent  (2.9) $\quad  \displaystyle 
\parallel \psi \parallel\ib{2} \,\,\, \leq \,\,\, G \, \parallel u \parallel\ib{0}
\,.$ 

\smallskip \noindent $\bullet \quad$ 
Then the field $\,\, \nabla \psi \,\,$ belongs to the space $\, H^1(\Omega) \, \times
\,  H^1(\Omega) .\, $  It is in consequence possible to interpolate this field in a
continuous way (see {\it e.g.}  Roberts and Thomas [RT91])  in the space $\, 
H({\rm div},\, \Omega) \,\,$ with the help of the fluxes on the edges : 

\smallskip \noindent  (2.10) $\quad  \displaystyle 
p_a \,\,= \,\, \int_{\displaystyle a} \, {{\partial \psi}\over{\partial n}} \,
{\rm d} \gamma \,,\qquad p\,=\, \sum_{a \, \in {\cal A}\ib{\cal T}} \, p_a \,\,
\varphi_a \, \,\, \in\,\,  H\ib{\cal T}({\rm div},\, \Omega)  \,$ 

\smallskip \noindent
and there exists a constant $\, L > 0 \,$ that only depends on the parameter $\,
\theta \,$ such that 

\smallskip \noindent  (2.11) $\quad  \displaystyle 
\parallel p \parallel \ib{H({\rm div},\, \Omega)} \,\,\, \leq \,\,\,L \,\, \parallel
u \parallel \ib{0} \,. \,$ 

\smallskip \noindent $\bullet \quad$ 
We observe that we have exactly 

\smallskip \noindent  (2.12) $\quad  \displaystyle 
{\rm div} \, p \,\, = \,\, u \qquad {\rm in} \,\, \Omega \,.\,$ 

\smallskip \noindent
On one hand, the two fields $\, {\rm div} \, p \,$ and $\, u \,$ are constant in each
element $\, K \,$ of the mesh $\, {\cal T} .\,$  On the other hand, we have :

\smallskip \noindent   $  \displaystyle 
\int_K \, {\rm div} \, p \,\,{\rm d}x \,\,= \,\, \int_{\partial K} \, p \smb n \,\,
{\rm d}\gamma \,\,=\,\,  \int_{\partial K} \,{{\partial \psi}\over{\partial n}} \,\,
{\rm d}\gamma \,\,=\,\, \int_K \,\Delta \psi \,\, {\rm d}x \,\,= \,\, \int_K \, u
\,\, {\rm d}x \,$ 

\smallskip \noindent
and the relation (2.12) is a consequence of the above property for the mean values. 

\smallskip \noindent $\bullet \quad$ 
Let $\, \, \Pi \, p \,$ be defined according to the Hypothesis 1, and 

\smallskip \noindent  (2.13) $\quad  \displaystyle 
q \,\,= \,\, {{1}\over{E}} \, \Pi \, p \,.\, $ 

\smallskip \noindent
We have as a consequence of (2.5) and (2.12) :  

\smallskip \noindent   $  \displaystyle 
( \, u \,,\, {\rm div} \, q \,) \,\,= \,\, {{1}\over{E}} \, ( \, {\rm div} \, p \,,\,
{\rm div} \,  \Pi \, p \,) \,\, \geq \,\, \parallel {\rm div} \, p \parallel\ib{0} ^2
\,\,\,= \,\,\, \parallel u \parallel\ib{0} ^2 \,$ 

\smallskip \noindent
that establishes (2.7). Moreover, we have due to (2.3), (2.4) and (2.11) : 

\smallskip \noindent   $  \displaystyle
\parallel q \parallel\ib{0} \,\,\,= \,\,\, {{1}\over{E}} \,\parallel   \Pi \, p
\parallel\ib{0} \,\,\, \leq \,\,\, {{B}\over{E}} \,\parallel    p \parallel\ib{0}
\,\,\, \leq \,\,\, {{B L }\over{E}} \,\parallel  u \parallel\ib{0} \, $

\smallskip \noindent   $  \displaystyle
\parallel {\rm div} \, q \parallel\ib{0} \,\,\, = \,\,\, {{1}\over{E}} \,\parallel
{\rm div} \, \Pi \, p \parallel\ib{0} \,\,\, \leq \,\,\, {{D}\over{E}} \,\parallel
{\rm div} \, p \parallel\ib{0}  \,\,\, = \,\,\,  {{D}\over{E}} \,\parallel  u
\parallel\ib{0} \,.\,  $

\smallskip \noindent 
Then due to the definition (1.6),  the two above inequalities establish 
the estimate (2.6) with $\,\, F \,=\,  {{1}\over{E}}  \,\sqrt{ B^2 L^2 + D^2} \,\,$ and
the Proposition is proven. $\hfill \square$ 

\bigskip \noindent
{\bf Proposition 2. \quad Discrete stability}

\noindent
Let $\, \theta \,$ be a strictly positive parameter. We suppose that the dual
Raviart-Thomas basis satisfies the Hypothesis 1. Then we have the following discrete
stability for the Petrov-Galerkin mixed formulation (1.36) : 

\setbox20=\hbox{ $ \displaystyle  
\exists \, \beta   > \,0 \,,\,\,\, \forall \, {\cal T} \in \, {\cal U}\ib{\theta}
\,,\,\, \forall \, \xi \, \in V\ib{\cal T} \,\, $ such that $ \displaystyle  \,\,
\,\, \parallel \xi \parallel\ib{V}  \,\,=\,\, 1 \,,\,\,$}
\setbox21=\hbox{ $ \displaystyle 
\exists \, \eta \, \in  \, V\ib{\cal T}^{\displaystyle \star}  \,,\,\, 
\parallel \eta \parallel\ib{V}  \,\,\leq\,\, 1 \,\,\, {\rm and} \,\,\, \gamma(\xi,\,
\eta) \,\, \geq \,\, \beta \,.\,$}
\setbox26= \vbox {\halign{# \cr \box20 \cr \box21 \cr   }}
\setbox27= \hbox{ $\vcenter {\box26} $}
\setbox28=\hbox{\noindent (2.14) $\quad \left\{  \box27 \right. $}
\smallskip \noindent $ \box28 $

\smallskip \noindent
with $\, \gamma (\smb,\, \smb) \,$ defined at the relation (1.11) and 
$\, \beta \,$ chosen such that 

\smallskip \noindent  (2.15) $\quad  \displaystyle 
\sqrt{ 1 - {{B + 2 D}\over{A}} \beta \,-\, \beta^2} \,\,\, \geq \,\,\, \biggl( \, 
1 + F \, \biggl( 1 \,+\, \sqrt{{{B + 2 A }\over{A}}} \, \biggr) \, \biggr) \,
\sqrt{\beta} \,. \,$ 

\smallskip \noindent
{\bf Proof of Proposition 2}. 

\smallskip \noindent $\bullet \quad$ 
We set $\,\,  \xi \equiv (u,\,p) \,\,$ satisfying the hypothesis (2.14) : 

\smallskip \noindent  (2.16) $\quad  \displaystyle 
 \parallel \xi   \parallel\ib{V}^2 \,\,\,\,\equiv \,\,\, \parallel u
\parallel\ib{0}^2 \,\,+\,\, \parallel p \parallel\ib{0}^2 \,\,+\,\, \parallel {\rm
div}\,p  \parallel\ib{0}^2  \,\,\, =\,\,\, 1  \,.\,$

\smallskip \noindent 
Then at last one of these terms is not too small and due to the three terms that arise
in relation (1.11), the proof is divided into three parts. 

\smallskip \noindent $\bullet \quad$ 
If the following condition 

\smallskip \noindent  (2.17) $\quad  \displaystyle  
\parallel {\rm div}\, p   \parallel\ib{0} \,\,  \geq \,\, \beta \,,\,$ 

\smallskip \noindent
is satisfied, we set  

\smallskip \noindent  (2.18) $\quad  \displaystyle  
v \,\,= \,\, {{{\rm div}\, p}\over{\parallel {\rm div}\, p  \parallel\ib{0}}} \, 
\,,\quad q \,\,=\,\, 0 \,\,,\qquad \zeta \,\,=\,\, (v,\,q) \, \in \,V\ib{\cal
T}^{\displaystyle \star} \,. \, $ 

\smallskip \noindent
Then $\,\, \, \displaystyle  \parallel {\rm div}\, v  \parallel\ib{0} \,=\, 1 \,\,
\,\,$ and   $\,\, \, \displaystyle  \parallel \zeta  \parallel\ib{0} \,\leq \, 1
\,.\, $ Moreover $\,\,\, \gamma(\xi,\, \zeta) \,=\, ( \, {\rm div} \, p \,,\, v \,)
\,\,= \,\, $ $\,\,= \,\, 
\parallel {\rm div}\, p  \parallel\ib{0} \,\, \geq \,\, \beta \,\,$ by
hypothesis (2.17) and the relation (2.14) is satisfied in this particular case. 

\smallskip \noindent $\bullet \quad$ 
Under  the following conditions 

\smallskip \noindent  (2.19) $\quad  \displaystyle  
\parallel {\rm div}\, p   \parallel\ib{0} \,\,  \leq \,\, \beta $ \quad and \quad 
$   \displaystyle   \parallel   p   \parallel\ib{0} \,\,\,\,  \geq \,\, \sqrt{{B + 2
D}\over{A}} \, \sqrt{\beta} \,\,, \,$

\smallskip \noindent
we set 

\smallskip \noindent  (2.20) $\quad  \displaystyle  
v \,\,= \,\, 0 \,\,\, \,,\quad q \,\,=\,\,  {{1}\over{B+D}} \,\, \Pi \, p 
 \,\,,\qquad \zeta \,\,=\,\, (v,\,q) \, \in \,V\ib{\cal
T}^{\displaystyle \star} \,. \, $ 

\smallskip \noindent
The following inequalities are a direct consequence of  (2.3) and (2.4) :  

\smallskip \noindent $ \displaystyle  
\parallel q  \parallel\ib{0} \,\,\, \leq \,\,\, {{B}\over{B+D}} \,\, 
\parallel p  \parallel\ib{0} \,$ \quad and \quad $ \displaystyle 
\parallel {\rm div}\, q   \parallel\ib{0} \,\,  \leq \,\, {{D}\over{B+D}} \,\,
\parallel {\rm div}\, p   \parallel\ib{0} \,$ 

\smallskip \noindent
so we deduce : 

\smallskip \noindent $ \displaystyle  
\parallel q \parallel  ^2 \ib{H({\rm div},\, \Omega)} \,\,\, \leq \,\,\, 
{{B^2 + D^2 }\over{(B+D)^2}} \,\, \parallel {\rm div}\, p   \parallel\ib{0} ^2 
\,\,\, \leq \,\,\, 1 \,$ 

\smallskip \noindent
because $\,\, B > 0 \,,\,\, D  > 0 \,. \,$  Then \quad $\, \parallel \zeta
\parallel\ib{V}  \,\,\leq\,\, 1 \,\, \,$  and we have also 

\smallskip \noindent $ \displaystyle  
\gamma (\xi,\, \zeta) \,\,= \,\, (p,\,q) \,+\, (u,\, {\rm div}
\, q) \,+\,  ({\rm div} \, p ,\,v ) \,$ 
\smallskip \noindent $ \displaystyle  \qquad \quad \,\,\,
\geq \,\,\, {{A}\over{B + D}} \, \parallel   p   \parallel\ib{0}^2 \,\,\,- \,\,\, 
 {{D}\over{B + D}} \, \, \beta \,\, \parallel u   \parallel\ib{0} \,$ 

\smallskip \noindent $ \displaystyle  \qquad \quad \,\,\,
\geq \,\,\,  {{1}\over{B + D}} \,\, \bigl( \, (B + 2 D) \, \beta \,-\, D \, \beta \,
\bigr) \,\, $ \qquad \qquad because $\,\,\,\, \parallel u   \parallel\ib{0} \,\leq \, 1
\,$ 

\smallskip \noindent $ \displaystyle  \qquad \quad \,\,\,
\geq \,\,\, \beta \,\,$ \qquad \qquad  \hfill 
and the property is established in this case. 

\smallskip \noindent $\bullet \quad$ 
If the two previous conditions (2.17) and (2.19) are in defect, {\it i.e.} if we have 

\smallskip \noindent  (2.21) $\quad  \displaystyle  
\parallel {\rm div}\, p   \parallel\ib{0} \,\,  \leq \,\, \beta $ \quad and \quad 
$   \displaystyle   \parallel   p   \parallel\ib{0} \,\,  \leq \,\, \sqrt{{B + 2
D}\over{A}} \, \sqrt{\beta} \,\,, \,$

\smallskip \noindent
then $\,\,\,  \displaystyle  \parallel u   \parallel\ib{0}^2 \,\,=\,\, 1 \,-\, 
\parallel p \parallel \ib{0} ^2 \,-\, \parallel {\rm div}\, p  \parallel \ib{0} ^2
\,\,\geq \,\, 1 - {{B + 2 D}\over{A}} \, \beta \,-\, \beta^2 \,\, \geq \,\, \beta \,>
\,0 \,$ 

\smallskip \noindent
due to the hypothesis (2.15). From the Proposition 1, there exists some vector field
$\,\, \widetilde{q} \in \, H\ib{\cal T}\sth ({\rm div},\, \Omega) \,\,$ satisfying
 $\quad  \displaystyle 
\parallel \widetilde{q} \parallel  \ib{H({\rm div},\, \Omega)} \,\,\, \leq \,\,\, F \,
\parallel u \parallel \ib{0} \,$ \quad and   $\quad  \displaystyle 
(\, u \,,\, {\rm div} \, \widetilde{q}  \, ) \quad \,  \geq \,\, \,\, \parallel u
\parallel \ib{0} ^2\, \,. \, $ \quad We set 

\smallskip \noindent  (2.22) $\quad  \displaystyle  
v \,\,= \,\, 0 \,\,\, \,,\quad q \,\,=\,\,  {{1}\over{F}} \,\, \widetilde{q} 
 \,\,,\qquad \zeta \,\,=\,\, (v,\,q) \, \in \,V\ib{\cal
T}^{\displaystyle \star} \,, \, $ 

\smallskip \noindent 
then 

\smallskip \noindent    $  \displaystyle  
\parallel \zeta \parallel\ib{V}  \,\,=\,\, {{1}\over{F}}\,
\parallel \widetilde{q} \parallel  \ib{H({\rm div},\, \Omega)} \,\,\, \leq\,\,\,  
\parallel u \parallel \ib{0}  \,\,\leq\,\, 1 \,$ 

\smallskip \noindent 
due to the hypothesis (2.14) relative to $\,   \parallel \xi \parallel\ib{V} \,. \,$
Moreover, we have 

\smallskip \noindent $ \displaystyle  
\gamma (\xi,\, \zeta) \,\,\,\,= \,\, (p,\,q) \,+\, (u,\, {\rm div}
\, q) \,+\,  ({\rm div} \, p ,\,v ) \,$ 
\smallskip \noindent $ \displaystyle  \qquad \qquad 
= \,\,\, {{1}\over{F}} \, (p,\, \widetilde{q} ) \,\,+\,\, {{1}\over{F}} \, (u ,\,
{\rm div} \, \widetilde{q} ) \,$ 
\smallskip \noindent $ \displaystyle  \qquad \qquad 
\geq \,\,\, {{1}\over{F}} \, \biggl( \, -\sqrt{{B + 2 D}\over{A}} \, \sqrt{\beta} \,\,
F \, \parallel u \parallel \ib{0} \,\,+\,\, \parallel u \parallel \ib{0} ^2 \,
\biggr) \,$ 
\smallskip \noindent $ \displaystyle  \qquad \qquad 
\geq \,\,\,\parallel u \parallel \ib{0} \,\, \biggl( \, {{1}\over{F}} \,\parallel u
\parallel \ib{0} \,-\, \sqrt{{B + 2 D}\over{A}} \, \sqrt{\beta} \, \biggr) \,$ 
\smallskip \noindent $ \displaystyle  \qquad \qquad 
\geq \,\,\, \sqrt{\beta} \,\, \biggl( \,  {{1}\over{F}} \, \biggl( \, 
1 + F \, \Bigl( 1 \,+\, \sqrt{{{B + 2 A }\over{A}}} \, \Bigr) \, \biggr) \,
\sqrt{\beta}\,- \, 
\sqrt{{B + 2 D}\over{A}} \, \sqrt{\beta} \, \biggr) \,$ 
\smallskip \noindent $ \displaystyle  \qquad \qquad 
\geq \,\,\, \sqrt{\beta} \,\, \Bigl( \,  {{1}\over{F}} \, \sqrt{\beta} \,+\,
\sqrt{\beta} \, \Bigr) \,$ 
\smallskip \noindent $ \displaystyle  \qquad \qquad 
\geq \,\,\,  \beta \,\, $ \qquad  \qquad \hfill 
and the property is satisfied for this last case.

\smallskip \noindent
The Proposition 2 is established.  $\hfill \square$ 

\bigskip \noindent
{\bf Theorem 1. \quad Error estimate}

\noindent
Let $\, \Omega \,$ be a two-dimensional open convex domain of $\, \R^2 \,$ with a
polygonal boundary,  $\, u \, \in H^2(\Omega) \,$ be the solution of the problem
(1.1)(1.2) considered under  variational formulation and $\, p = \nabla u \,$ be the
associated momentum. Let $\, \theta \,$ be a strictly positive parameter,   $\,
{\cal U}\ib{\theta} \,$  a family of meshes $\, {\cal T} \,$ and $\,  V\ib{\cal
T}^{\displaystyle \star} \,$ defined in (1.35) and associated with a choice of a dual
Raviart-Thomas basis that satisfies the  Hypothesis 1. Let $\, \xi\equiv (u\ib{\cal
T},\,p\ib{\cal T}) \, \in  V\ib{\cal T} \,$ be the solution of the discrete problem
(1.34). Then there exists some constant  $\, C > 0 \,$ that only depends on the
parameter $\, \theta \,$ such that 

\smallskip \noindent  (2.23) $\quad  \displaystyle  
\parallel u - u \ib{\cal T} \parallel \ib{0} \,+\, \parallel p - p \ib{\cal T}
\parallel  \ib{H({\rm div},\, \Omega)} \,\,\,\,   \leq  \, \,\, C \, \, h\ib{\cal T}
\,\,  \parallel f \parallel \ib{0} \,\,.  \,$ 

\smallskip \noindent
{\bf Proof of Theorem 1}. 

\smallskip \noindent $\bullet \quad$ 
On one hand, it is sufficient to apply the general approximation Theorem established
by Babu\u ska's  for continuous (respectively discrete) variational mixed systems 
(1.13) (respectively (1.36)) {\it i.e.} to verify that the bilinear form $\,
\gamma(\smb,\, \smb)\,$ defined in (1.11) is continuous on the Hilbert space $\, V
\,=\, L^2(\Omega) \,\times\,  H({\rm div},\, \Omega) \,,\,$ which is clear. It is also 
necessary  to verify the so-called discrete  inf-sup condition (2.14), that has been
established at the Proposition 2. Last but not least, it is necessary to satisfy the
following infinity condition :

\smallskip \noindent  (2.24) $\quad  \displaystyle 
\forall \, \eta \in V\ib{\cal T}^{\displaystyle \star} \,,\quad \eta \ne 0 \,\,
\Longrightarrow \,\, \sup_{\xi \in  V\ib{\cal T}}\, \gamma (\xi,\,\eta) \,=\, +\infty
\,. \,$ 

\smallskip \noindent $\bullet \quad$ 
The infinity condition (2.24) is established as follows. Let $\, \zeta \equiv
(v,\,q) \in   V\ib{\cal T}^{\displaystyle \star} \,$ be a ``test vector''  different
from zero. If there exists some mesh element $\, K \in {\cal E}\ib{\cal T} \,$ such
that $\, \int_{K} {\rm div} \,q \, {\rm d}x \, \ne \, 0,\,$ then we consider $\, \xi =
(u,\,p) \,$ chosen according to $\, u = \lambda \, \widetilde{u} \,$ and $\, p = 0 .\,$
We suppose that the  field $\,  \widetilde{u} \in L^2\ib{\cal T} \,$ is null for all
the elements of the mesh $\, {\cal T} \,$ except for the particular element $\, K \,$
where we suppose $\, \widetilde{u}\ib{K} = \int_{K} {\rm div} \,q \, {\rm d}x .\,$ Then
we have  $\,\, \, \gamma(\xi,\, \zeta) \equiv (p,\,q) + (u,\, {\rm div}\,q) \,+\, ({\rm
div}\, p,\,v) \,=\, \lambda \, \smash{ \bigl( \int_{K} {\rm div} \,q \, {\rm d}x
\bigr)^2 } ,\,\, \,$ which tends to infinity as $\, \lambda \,$ tends to infinity. If  $\, \int_{K} {\rm div}
\,q \, {\rm d}x = 0 \,$ for all mesh elements $\, K \in {\cal E}\ib{\cal T} \,$ 
and if the field $\, q \,$ is not null, we can write it on the form $\,\, q \,=\, \Pi
\, \widetilde{p} \,\,$ with $\,\,  \widetilde{p} \, \in  H\ib{\cal T}({\rm div},\,
\Omega) \,\,$ because the mapping $\,\, \Pi \,\,$ is clearly bijective due to the
property (2.2).   We set  $\,\,\,  p = \lambda \, \widetilde{p} \,\, \,$  and 
$\, u = 0 .\,$ Then $\, \gamma(\xi,\zeta) \,=\,(p,\,q) \,=\, 
\lambda \, (\widetilde{p} ,\,q) \,\, \geq \,\, \lambda \, A \, \parallel q
\parallel\ib{0}^2 \, \,$ due to the hypothesis (2.2)~;   the infinity property
(2.24) is established in this second particular case because $\, q \ne 0 .\,$ 
 If $\, q = 0 ,\,$ then $\, v \,$ is not null due to the left hand side of (2.24).
Following the proof of Proposition 1, we introduce the vector field $\, \widetilde{p} 
\,\in H\ib{\cal T}({\rm div},\, \Omega) \,\,$   satisfying  the relations (2.11) and
(2.12) :  $\,\,\,  {\rm div} \,\widetilde{p} \equiv v \,\,$ and $\,\,  \parallel
\widetilde{p} \parallel \ib{H({\rm div},\, \Omega)} \,  \leq  \,L  \,  \parallel v 
\parallel \ib{0} . \,$  We set  $\,\,\,  p = \lambda \, \widetilde{p} ,\,$   
$\, u = 0 \,$ and $\, \xi = (u,\,p).\,$ 
Then  $\, \gamma(\xi,\zeta) \,=\,\lambda \, ( {\rm div} \, \widetilde{p}
,\, v) \,$ $\,=\, \lambda \, \parallel v \parallel\ib{0}^2 \, \,$ tends to infinity
when $\, \lambda \,$ tends to infinity, and the infinity condition (2.24) is
established. 

\smallskip \noindent $\bullet \quad$ 
The conclusion of the Babu\u ska's Theorem [Ba71] assures the existence of
some constant $\, C > 0 \,$ that only depends on $\, \theta \,$ such that the error
between the solution of the continuous problem (1.13) and the discrete problem (1.23)
is majorated by the interpolation error : 

\setbox20=\hbox{ $ \displaystyle  
\parallel u - u \ib{\cal T} \parallel \ib{0} \,+\, \parallel p - p \ib{\cal T}
\parallel  \ib{H({\rm div},\, \Omega)} \,\,\,\,   \leq  \, \,\, $}
\setbox21=\hbox{ $ \displaystyle 
\leq  \,\, C \, \,\, \Bigl( \, \inf_{ \displaystyle   v \in  L^2\ib{\cal T}} 
\parallel u - v \parallel \ib{0} \,\,+\,\, \inf_{ \displaystyle   q
\in  H\ib{\cal T}({\rm div},\, \Omega) }  \parallel p - q \parallel  \ib{H({\rm div},\,
\Omega)}  \, \Bigr)  \,. \,$}
\setbox26= \vbox {\halign{# \cr \box20 \cr \box21 \cr   }}
\setbox27= \hbox{ $\vcenter {\box26} $}
\setbox28=\hbox{\noindent (2.25) $\quad \left\{  \box27 \right. $}
\smallskip \noindent $ \box28 $

\smallskip \noindent
Then following classical interpolation results for scalar [CR72] and vectorial [RT77]
fields, we deduce from (2.25) :

\smallskip \noindent   $  \displaystyle 
\parallel u - u \ib{\cal T} \parallel \ib{0} \,+\, \parallel p - p \ib{\cal T}
\parallel  \ib{H({\rm div},\, \Omega)} \,\,\,\,   \leq  \, \,\, C \, \,\, \, \bigl( \, 
h\ib{\cal T} \,\parallel u  \parallel \ib{1} \,\,+\,\, 
 h\ib{\cal T} \,\parallel p  \parallel \ib{1} \, \bigr)  \, \,$
 
\smallskip \noindent   $\qquad  \qquad  \qquad   \displaystyle 
\,\,\,\,   \leq  \, \,\, C \quad h\ib{\cal T}  \parallel u 
\parallel \ib{2} \qquad \leq  \, \,\, C \quad  h\ib{\cal T} \, \parallel f
\parallel \ib{0}  \, \,$

\smallskip \noindent 
and the Theorem 1 is proven.  $\hfill \square$


\bigskip
\bigskip
\noindent {\smcaps 3)$\,\,\,$ Towards a first Petrov-Galerkin finite volume scheme } 

\smallskip \noindent $\bullet \quad$ 
We propose in this section to formulate some ideas in order to construct  a dual
Raviart-Thomas basis  $\, \varphi^{\star}\ib{a} \, $ where $\,\, a \,\,$ is 
an internal edge of the mesh $\, {\cal T }\,$ $( a \in {\cal A}\ib{\cal T} ) .\,$  
Following (1.15) and (1.16),  we denote by $\, a \equiv ({\rm S},\, {\rm N}) \,$ this
edge, by $\, {\rm O} \,$  the middle of $\, {\rm SN}  \,$ 
and by $\, K,\,L \,$ the two triangles that compose the  co-boundary. The normal
$\, n_a \,$ is supposed to be oriented from the element $\, K \,$ towards the element
$\, L \,$ and there exists two vertices $\, {\rm W }\,$ and  $\, {\rm
E }\,$ such that the relation (1.17) holds (see the Figure 1). We consider the four
edges  $\, ({\rm N},\, {\rm W}), \,$  $\, ({\rm W},\, {\rm S}), \,$
 $\, ({\rm S},\, {\rm E}) \,$ and $\, ({\rm E},\, {\rm N})  \,$   that compose the
boundary of the union $\, K \cup L .\,$ We define four new triangles $\, M,\, P, \, Q
\, $ and $\, R \,$ and four new vertices $\,\, {\rm A},\, {\rm B},\, {\rm C}\, $ and
$ \, {\rm D}\, \,$   in the mesh $\, {\cal T} \,$  by the relations

\setbox20=\hbox{ $ \displaystyle  
\partial^c({\rm E},\, {\rm N}) \,\,\,  \equiv (L,\,M) \,\,, \qquad M \, \equiv \, 
({\rm N},\,{\rm E},\,{\rm A}) \,  $}
\setbox21=\hbox{ $ \displaystyle 
\partial^c({\rm N},\, {\rm W}) \, \equiv (K,\,P) \,\,\,\, ,\qquad P \, \equiv \, 
({\rm W},\,{\rm N},\,{\rm B}) \,  $}
\setbox22=\hbox{ $ \displaystyle  
\partial^c({\rm W},\, {\rm S}) \,\, \equiv (K,\,Q) \,\,\,\,,\qquad Q \, \equiv \, 
({\rm S},\,{\rm W},\,{\rm C}) \,  $}
\setbox23=\hbox{ $ \displaystyle  
\partial^c({\rm S},\, {\rm E}) \,\,\,\, \equiv (L,\,R) \,\,\,\,\,\, ,\qquad R \,
\equiv \,  ({\rm E},\,{\rm S},\,{\rm D}) \,  $}
\setbox26= \vbox {\halign{# \cr \box20 \cr \box21 \cr \box22 \cr \box23 \cr  }}
\setbox27= \hbox{ $\vcenter {\box26} $}
\setbox28=\hbox{\noindent (3.1) $\quad \left\{  \box27 \right. $}
\smallskip \noindent $ \box28 $

\smallskip \noindent
as illustrated on the Figure 2.

\bigskip 
\centerline { \epsfysize=5,0cm    \epsfbox  {Petrov.epsf} }
\smallskip  \smallskip 

\centerline {    {\bf Figure 2} : support $\, {\cal V}({\rm S},\,{\rm N}) \,$ of the
dual Raviart-Thomas basis function $\, \varphi^{\star}\ib{\rm SN} .\,$  }
\smallskip 

\bigskip 
\vfill \eject
\noindent
{\bf Hypothesis 2. \quad A simple choice for  dual Raviart Thomas basis functions}

\noindent
We suppose first that the Raviart Thomas dual basis $\, \varphi^{\star}\ib{b} ,\,$
$\, ( b \in  {\cal A}\ib{\cal T} ) \,$  satisfies the $\, H({\rm div}) $ conformity
property  (1.31) and the orthogonality   (1.32). Moreover, we suppose that for  each
internal edge $\, a \equiv ({\rm S},\, {\rm N}) ,\,$  the  support of the dual
Raviart-Thomas basis function $\, \varphi^{\star}\ib{\rm SN} \,$ is included in a
vicinity $\, {\cal V}(a) =  {\cal V}({\rm S},\, {\rm N})  \,$  composed by  the six
triangles  $\,\,  K ,\,  L  ,\, M  ,\, P  ,\,   Q  \,$ and $ \,\,   R    \,\,$ 
introduced previously (see the Figure 2) : 

\setbox20=\hbox{ $ \displaystyle  
 {\cal V}({\rm S},\, {\rm N}) \,\,\subset  \,\,  K \,\cup \, L \,\cup \, M \,\cup \,P
  \,\cup \,  Q  \,\cup \,   R    \,,\,\,$ }
\setbox21=\hbox{ $ \displaystyle 
{\rm supp} \, \bigl( \,  \varphi^{\star}\ib{\rm SN}  \, \bigr) \,\, \subset \,\, 
 {\cal V}({\rm S},\, {\rm N}) \,.\   $}
\setbox26= \vbox {\halign{# \cr \box20 \cr \box21 \cr   }}
\setbox27= \hbox{ $\vcenter {\box26} $}
\setbox28=\hbox{\noindent (3.2) $\quad \left\{  \box27 \right. $}
\smallskip \noindent $ \box28 $
 
\smallskip \noindent
We suppose also that the divergence field  $\,\,\, {\rm div} \, \varphi^{\star}\ib{a}
\,\,\,$ is {\bf constant in each triangle}   of the mesh : 

\smallskip \noindent  (3.3) $\quad  \displaystyle 
{\rm div} \, \varphi^{\star}\ib{a} \,\, \in \,\, L\ib{\cal T}^2(\Omega) \,\,,\qquad 
\forall \, a \, \in  {\cal A}\ib{\cal T} \,.  \,$ 

\bigskip 
\centerline { \epsfysize=6cm    \epsfbox  {fleches.epsf} }
\smallskip  \smallskip 
\centerline {    {\bf Figure 3} : Notations and orientations.  }
\smallskip

\bigskip \noindent
{\bf Theorem 2. \quad Necessary condition for a dual Raviart-Thomas basis. }

\noindent
Let  $\,\, \varphi^{\star}\ib{\rm SN}  \,\,$ be a dual Raviart Thomas basis function
satisfying the Hypothesis 2. We introduce the following fluxes accross the internal
edges $\, {\rm S}{\rm N} ,\,$ $\, {\rm E}{\rm N} ,\,$ $\, {\rm N}{\rm W} ,\,$
$\, {\rm W}{\rm S} \,$ and  $\,{\rm S}{\rm E} \,$ respectively : 
\setbox20=\hbox{ $ \displaystyle   \eta  \, \equiv \,  \int_{\rm SN} \,
\varphi^{\star}\ib{\rm SN}  \, \smb \, n\ib{\rm SN} \,\,\, {\rm d}\gamma  \,,\,\,$ }
\setbox21=\hbox{ $ \displaystyle   \alpha \, \equiv \,  \int_{\rm EN} \,
\varphi^{\star}\ib{\rm SN}  \, \smb \, n\ib{\rm EN}  \,\,\,  {\rm d}\gamma \,,\,\,$ 
 $ \displaystyle \quad   \beta \, \equiv \,  \int_{\rm NW}
\,\varphi^{\star}\ib{\rm SN}  \, \smb \, n\ib{\rm NW}  \,\,\, {\rm d}\gamma \,,\,\,$ }
\setbox23=\hbox{ $ \displaystyle   \gamma \, \equiv \,  \int_{\rm WS} \, 
\varphi^{\star}\ib{\rm SN}  \, \smb \, n\ib{\rm WS}  \,\,\,  {\rm d}\gamma  \,,\,\,$ 
 $ \displaystyle \quad   \delta \, \equiv \,  \int_{\rm SE} \,
\varphi^{\star}\ib{\rm SN}  \, \smb \, n\ib{\rm SE}  \,\,\, {\rm d}\gamma \, .\,\,$ }
\setbox26= \vbox {\halign{# \cr \box20 \cr \box21 \cr \box23 \cr }}
\setbox27= \hbox{ $\vcenter {\box26} $}
\setbox28=\hbox{\noindent (3.4) $\quad \left\{  \box27 \right. $}
\smallskip \noindent $ \box28 $
\smallskip \noindent 
Then we have the necessary conditions : 
\smallskip \noindent  (3.5) $\quad  \displaystyle
\eta  \,   \fl {\rm KL}    \,+\, 
\alpha  \, {\fl {\rm LM}}    \,+\, 
\beta  \,  {\fl {\rm KP}}    \,+\, 
\gamma  \, {\fl {\rm KQ}}    \,+\,
\delta  \, {\fl {\rm LR}}   \,\,=\,\, \abs{ \fl {\rm SN}}  \,  n\ib{\rm SN}  \,$
\smallskip 
\setbox20=\hbox{ $ \displaystyle   
\alpha  \,  \fl {\rm LM}  \,\smb\,   \fl {\rm WA}  \,+\, 
\beta   \,  \fl {\rm KP}  \,\smb\,   \fl {\rm EB}  \,+\,   
\gamma  \,  \fl {\rm KQ}  \,\smb\,   \fl {\rm EC}  \,+\, 
\delta  \,  \fl {\rm LR}  \,\smb\,   \fl {\rm WD}  \,\,= \,\,   $ }
\setbox21=\hbox{ $ \displaystyle  \qquad \qquad   \qquad \qquad   \qquad \qquad  
\quad   \,\,=\,\,  -3   \abs{  \fl{\rm SN}} \,    n\ib{\rm SN}  \,\smb\,   \bigl(  
\fl {\rm OL} + \fl{\rm OK}  \bigr) \,  . \, $ }
\setbox26= \vbox {\halign{# \cr \box20 \cr \box21 \cr  }}
\setbox27= \hbox{ $\vcenter {\box26} $}
\setbox28=\hbox{\noindent (3.6) $\quad \left\{  \box27 \right. $}
\smallskip \noindent $ \box28 $

\smallskip \smallskip \noindent $\bullet \quad$ 
The finite volume approach  is then obained with a six point scheme for the mean
gradient in the normal direction in the manner of (1.29)  thanks to the first equation
of the mixed variational formulation (1.24) : 

\setbox20=\hbox{ $ \displaystyle   
\int_{\rm SN}    \nabla  u\ib{\cal T} \, \smb n \,\, \, {\rm d}\gamma \,\,= \,\,     
\eta \, (u\ib{\rm L} - u\ib{\rm K}) \,+\, $ }
\setbox21=\hbox{ $ \displaystyle  \quad + \,  
\alpha \, (u\ib{\rm M} -  u\ib{\rm L}) \,+\,  \beta \, (u\ib{\rm P} -  u\ib{\rm K}) 
\,+\,  \gamma \, (u\ib{\rm Q} -  u\ib{\rm K}) \,+\,   \delta \, (u\ib{\rm R} - 
u\ib{\rm L}) .\,$ }
\setbox26= \vbox {\halign{# \cr \box20 \cr \box21 \cr  }}
\setbox27= \hbox{ $\vcenter {\box26} $}
\setbox28=\hbox{\noindent (3.7) $\quad \left\{  \box27 \right. $}
\smallskip \noindent $ \box28 $

\smallskip \noindent 
We remark that the constraints (3.5) express that the relation (3.7) is {\bf exact} if
the field $\, u\ib{\cal T} \,$ is an affine function.  

\bigskip 
\centerline { \epsfysize=3,5cm    \epsfbox  {edge.epsf} }
\smallskip  \smallskip
 
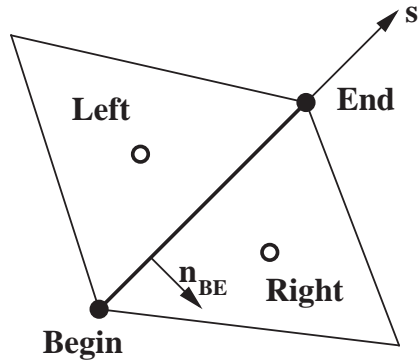
\centerline {    {\bf Figure 4} : Notations for an arbitrary edge {\rm B}{\rm E}.  }
\smallskip 

\smallskip \noindent $\bullet \quad$ 
We precise some notations that we will use in the next pages. Let $\, ({\rm B},\,{\rm
E})\,$ be an edge of the mesh (see {\it e.g.} the Figure 4), and $\, (L,\,R )\,$ its
co-boundary. If the edge is directed from $\, {\rm B} \,$ towards  $\, {\rm E} ,\,$
the axis $\, s\,$ has its origin at vertex  $\, {\rm B} \,$ and the normal $\,
n\ib{\rm BE} \,$ is oriented from $\, L \,$ to $\, R \,$ in such a way that the pair of
vectors $\,\, ( n\ib{\rm BE} ,\, \fl{\rm BE})\,$ is direct. If $\,\, \xi \equiv 
\int_{\rm BE} \, \varphi^{\star}  \, \smb \, n\ib{\rm BE} \,\,\, {\rm d}s \,\,$
is the flux of the function $\, \varphi^{\star} \,$ accross the edge $\, ({\rm
B},\,{\rm E}),\,$ we will denote by $\, \xi\ib{1} \,,\,$ $ \,
\smash{\widetilde{\xi}}\ib{1}  \,,\,$ $\, \xi\ib{2} \, ,\,$ and $\,
\smash{\widetilde{\xi}}\ib{2} \, \,$ the following momenta : 

\setbox20=\hbox{ $ \displaystyle   
 \xi\ib{1} \,=\, \int_{\rm BE} \, \, \varphi^{\star}   \smb \, n\ib{\rm BE}
\,\,\,s \, {\rm d}s \,,\quad 
\smash{\widetilde{\xi}}\ib{1} \,=\, \int_{\rm BE} \, \, \varphi^{\star}   \smb
\, n\ib{\rm BE} \,\,\,({\rm BE} - s)  \, {\rm d}s \,,\,\,  $ }
\setbox21=\hbox{ $ \displaystyle  
 \xi\ib{2} \,=\, \int_{\rm BE} \, \, \varphi^{\star}   \smb \, n\ib{\rm BE}
\,\,\,s^2 \, {\rm d}s \,,\quad 
\smash{\widetilde{\xi}}\ib{2}\,=\, \int_{\rm BE} \, \, \varphi^{\star}   \smb \,
n\ib{\rm BE} \,\,\,({\rm BE} - s)^2  \, {\rm d}s \,.\,\,  $ }
\setbox26= \vbox {\halign{# \cr \box20 \cr \box21 \cr  }}
\setbox27= \hbox{ $\vcenter {\box26} $}
\setbox28=\hbox{\noindent (3.8) $\quad \left\{  \box27 \right. $}
\smallskip \noindent $ \box28 $

\smallskip \smallskip \smallskip \noindent 
The proof of Theorem 2 needs a certain number of technical lemmae and preliminary
propositions.

\bigskip  \noindent
{\bf Lemma 1. \quad Radius of giration } 

\noindent
Let $\, M = ({\rm N},\, {\rm E},\,{\rm A}) \,$ be a triangle of the mesh $\, {\cal
T} \,$ and $\, {\rm M}\,$ its associated center of gravity (see the Figure 4). We will
denote by $\, \rho\ib{M}\,$ the radius of giration :
\smallskip \noindent  (3.9) $\quad  \displaystyle
\rho\ib{\! M} \,\,= \,\, \sqrt{ \, 
{{1}\over{36}} \, ( {\rm AN}^2 \,+\, {\rm NE}^2 \,+\, {\rm EA}^2  \,) \,}  \,$  
\smallskip \noindent
Then we have 
\smallskip \noindent  (3.10) $\quad  \displaystyle
{{1}\over{\abs{M}}} \, \int_{M} \, \abs{x-{\rm N}}^2 \, {\rm d} x \,\, = \,\,
\rho\ib{\! M}^2 \, +\,  {\rm NM}^2 \,.\,$

$~$
 
\bigskip 
\centerline { \epsfysize=3,5cm    \epsfbox  {triangle.epsf} }
\smallskip  \smallskip 

\centerline {    {\bf Figure 5} : About the radius of giration.  }
\smallskip \smallskip \smallskip

\smallskip \noindent
{\bf Proof of Lemma 1}. 

\smallskip \noindent 
We have on one hand : 

\smallskip \noindent  $ \displaystyle
\abs{M} \,=\, {{1}\over{2}} \,h \, (p+q) \,,\, \,$  

\smallskip \noindent  $ \displaystyle
\int_{M} \, \abs{x-{\rm N}}^2 \, {\rm d} x \,\, = \,\, \int_{\displaystyle
0}^{\displaystyle h} \,{\rm d}y \,  \,   \int_{\displaystyle -p \, {{y}\over{h}}
}^{\displaystyle q{{y}\over{h}}} \, \,{\rm d} x \, \, \bigl( x^2 + y^2 \bigr)  \,$ 
\smallskip \noindent  $ \displaystyle \qquad \qquad \qquad  \quad  \,\,\,\,\, = \,\, 
\int_{\displaystyle 0}^{\displaystyle h} \,{\rm d}y \, \biggl( \, {{1}\over{3}} \,
(p^3 + q^3) \, \Bigl({{y}\over{h}} \Bigr)^3 \, +\, (p+q) \, {{y}\over{h}} \, y^2 \,
\biggr) \,  $
\smallskip \noindent  $ \displaystyle \qquad \qquad \qquad  \quad  \,\,\,\,\, = \,\, 
{{1}\over{12}} \, (p+q) \, (p^2 - p\,q + q^2) \, h \,+\, {{1}\over{4}} \, (p+q) \, h^3
\,$ 
\smallskip \noindent  $ \displaystyle \qquad \qquad \qquad  \quad  \,\,\,\,\, = \,\, 
\abs{M} \, \Bigl( \, {{1}\over{6}} \, (p^2 - p\,q + q^2) \,+\, {{1}\over{2}} \, h^2
\, \Bigr) \,$ 

\smallskip \noindent 
and on the other hand : 

\smallskip \noindent  $ \displaystyle
{{1}\over{36}} \,  ( {\rm AN}^2 \,+\, {\rm NE}^2 \,+\, {\rm EA}^2 \,) \,+\, {\rm
NM}^2 \,\,=\,\, $ 
\smallskip \noindent  $ \displaystyle \qquad   \quad  = \,\, 
{{1}\over{36}} \, \bigl[ \, p^2 + h^2 + q^2 + h^2 + (p+q)^2 \, \bigr]
\,+\, {{1}\over{9}} \, \bigl[ \,(p-q)^2 \,+\, 4 \,h^2 \, \bigr]  \,$ 
\smallskip \noindent  $ \displaystyle \qquad   \quad  = \,\, 
{{1}\over{36}} \, \bigl( \, 6 \, p^2 \,+\, 6 \, q^2 \,-\, 6 \, p \, q \, +\, 18 \,
h^2 \, \bigr) \,$ 
\smallskip \noindent  $ \displaystyle \qquad   \quad  = \,\, 
{{1}\over{6}} \, \bigl( \, p^2 \,+\, q^2 \,-\,p \, q \, +\,3 \, h^2 \, \bigr)\,.  \,$ 

\smallskip \noindent 
So the relation (3.10) is established. $\hfill \square$ 

\bigskip \noindent
{\bf Proposition 3.  \quad First relations between momenta }

\noindent
The Hypothesis 2 implies the following relations inside the triangle $\, M = ({\rm
N},\, {\rm E},\,{\rm A}) \,: \,$
\setbox20=\hbox{ $ \displaystyle   
\alpha\ib{1} \,=\,   \fl{\rm EM} \,\smb\, {{\fl{\rm EN}}\over{\rm
EN}} \,\, \alpha \,\,, \qquad \quad  \,\,\,\,\, 
\widetilde{\alpha}\ib{1}   \,=\, - \fl{\rm NM} \,\smb\, {{\fl{\rm
EN}}\over{\rm EN}} \,\, \alpha \, ,\,$  }
\setbox21=\hbox{ $ \displaystyle   
\alpha\ib{2} \,=\,  ( \,\rho\ib{\! M}^2 \, +\,  {\rm EM}^2 \,) \, \alpha \,\,, \qquad 
\smash{\widetilde{\alpha}}\ib{2} \,=\, ( \,\rho\ib{\! M}^2 \, +\,  {\rm NM}^2 \,) \,
\alpha \,, \,$  }
\setbox26= \vbox {\halign{# \cr \box20 \cr \box21 \cr  }}
\setbox27= \hbox{ $\vcenter {\box26} $}
\setbox28=\hbox{\noindent (3.11) $\quad \left\{  \box27 \right. $}
\smallskip \noindent $ \box28 $
\smallskip \noindent
and the analogous ones obtained from the Figure 3 in the triangles  
$\, P = ({\rm W},\, {\rm N},\,{\rm B}) ,\,  \,$
$\, Q = ({\rm S},\, {\rm W},\,{\rm C})  \,  \,$ and 
$\, R = ({\rm E},\, {\rm S},\,{\rm D}) \, : \,$
\setbox20=\hbox{ $ \displaystyle   
\beta\ib{1} \,=\,   \fl{\rm NP} \,\smb\, {{\fl{\rm NW}}\over{\rm
NW}} \,\, \beta \,\,, \qquad \quad \,\,
\widetilde{\beta}\ib{1}   \,=\, - \fl{\rm WP} \,\smb\, {{\fl{\rm
NW}}\over{\rm NW}} \,\, \beta \, , \,$  }
\setbox21=\hbox{ $ \displaystyle   
\beta\ib{2} \,=\,  ( \,\rho\ib{\! P}^2 \, +\,  {\rm NP}^2 \,) \, \beta \,\,, \qquad 
{\widetilde{\beta}}\ib{2} \,=\, ( \,\rho\ib{\! P}^2 \, +\,  {\rm WP}^2 \,) \,
\beta \,, \,$  }
\setbox26= \vbox {\halign{# \cr \box20 \cr \box21 \cr  }}
\setbox27= \hbox{ $\vcenter {\box26} $}
\setbox28=\hbox{\noindent (3.12) $\quad \left\{  \box27 \right. $}
\smallskip \smallskip \smallskip \noindent $ \box28 $
\setbox20=\hbox{ $ \displaystyle   
\gamma\ib{1} \,=\,   \fl{\rm WQ} \,\smb\, {{\fl{\rm WS}}\over{\rm
WS}} \,\, \gamma \,\,, \qquad \quad \,\,\,
\widetilde{\gamma}\ib{1}   \,=\, -  \fl{\rm SQ} \,\smb\, {{\fl{\rm
WS}}\over{\rm WS}} \,\, \gamma \, , \,$  }
\setbox21=\hbox{ $ \displaystyle   
\gamma\ib{2} \,=\,  ( \,\rho\ib{\! Q}^2 \, +\,  {\rm WQ}^2 \,) \, \gamma \,\,, \qquad 
\smash{\widetilde{\gamma}}\ib{2} \,=\, ( \,\rho\ib{\! Q}^2 \, +\,  {\rm SQ}^2 \,) \,
\gamma \,, \,$  }
\setbox26= \vbox {\halign{# \cr \box20 \cr \box21 \cr  }}
\setbox27= \hbox{ $\vcenter {\box26} $}
\setbox28=\hbox{\noindent (3.13) $\quad \left\{  \box27 \right. $}
\smallskip \smallskip \smallskip \noindent $ \box28 $
\setbox20=\hbox{ $ \displaystyle   
\delta\ib{1} \,=\,   \fl{\rm SR} \,\smb\, {{\fl{\rm SE}}\over{\rm
SE}} \,\, \delta \,\,, \qquad \qquad 
\widetilde{\delta}\ib{1}   \,=\, -  \fl{\rm ER} \,\smb\, {{\fl{\rm
SE}}\over{\rm SE}} \,\, \delta \, , \,$  }
\setbox21=\hbox{ $ \displaystyle   
\delta\ib{2} \,=\,  ( \,\rho\ib{\! R}^2 \, +\,  {\rm SR}^2 \,) \, \delta \,\,, \qquad
\,\,  \smash{\widetilde{\delta}}\ib{2} \,=\, ( \,\rho\ib{\! R}^2 \, +\,  {\rm ER}^2 \,)
\, \delta \,. \,$  }
\setbox26= \vbox {\halign{# \cr \box20 \cr \box21 \cr  }}
\setbox27= \hbox{ $\vcenter {\box26} $}
\setbox28=\hbox{\noindent (3.14) $\quad \left\{  \box27 \right. $}
\smallskip \smallskip \smallskip \noindent $ \box28 $

\smallskip \noindent
{\bf Proof of Proposition  3}. 

\smallskip \noindent $\bullet \quad$ 
We write the orthogonality (1.32) between the two edges $\, a =  ({\rm S},\,{\rm N}) 
\,$ and  the edge $\, b = ({\rm A},\,{\rm N}) \,$  (see the Figure 3). 
Inside the triangle  $\, M = ({\rm N},\, {\rm E},\,{\rm A}) ,\,$ we have 

\smallskip \noindent  $   \displaystyle
\varphi\ib{\rm AN}  \,\,=\,\,  {{1}\over{2 \abs{M}}} \, (x-{\rm E}) \,\,=\,\,
{{1}\over{4 \abs{M}}} \,\, \nabla (\abs{ x-{\rm E}}^2) \,$

\smallskip \noindent 
then 

\smallskip \noindent  $   \displaystyle
0 \,\,=\,\, \int_{\Omega} \, \varphi^{\star}\ib{\rm SN} \,\smb\,  \varphi\ib{\rm AN} \,
{\rm d}x   \,\,=\,\,  \int_{M} \, \varphi^{\star}\ib{\rm SN} \,\smb\,  \varphi\ib{\rm
AN} \, {\rm d}x  \,   $ 
\smallskip \noindent  $   \displaystyle     = \,\, 
- \, \int_{M} \, \bigl( {\rm div} \, \varphi^{\star}\ib{\rm SN} \bigr) \,\, 
{{1}\over{4 \abs{M}}} \,\,\abs{ x-{\rm E}}^2 \, {\rm d}x   \,\,\,+\,\,\, 
 \int_{\partial M} \, \bigl(   \varphi^{\star}\ib{\rm SN} \,\smb\,n \bigr) \,\,
{{1}\over{4 \abs{M}}} \,\,\abs{ x-{\rm E}}^2 \, {\rm d}\gamma \,$ 
\smallskip \noindent  $   \displaystyle     = \,\, 
-   \bigl( {\rm div} \, \varphi^{\star}\ib{\rm SN} \bigr)(M) \,  \int_{M}
{{1}\over{4 \abs{M}}} \,  \abs{ x-{\rm E}}^2 \, \, {\rm d}x    \,$ 
\smallskip \noindent  $   \displaystyle     \qquad \qquad  \qquad \qquad    - \,\, 
\,{{1}\over{4 \abs{M}}} \, 
  \int_{\rm NE}  \,  \bigl(   \varphi^{\star}\ib{\rm SN} \,\smb\,n\ib{\rm NE}
\bigr)  \,  \abs{ x-{\rm E}}^2  \, \, {\rm d}\gamma \,$ 
\smallskip \noindent  $   \displaystyle     = \,\, {{1}\over{4 \abs{M}}} \biggl( \, 
 {{\alpha}\over{\abs{M}}} \, \, \int_{M} \,\abs{ x-{\rm E}}^2 \, {\rm d}x   \,\,-\,\, 
\int_{\rm E N}   \bigl(   \varphi^{\star}\ib{\rm SN} \,\smb\,n\ib{\rm NE} \bigr) 
\,\abs{ x-{\rm E}}^2 \, {\rm d}s \, \, \biggr) \,$ 
\smallskip \noindent  $   \displaystyle     = \,\, {{1}\over{4 \abs{M}}} \, \Bigl( \, 
\alpha \,\,  (\, \rho\ib{\! M}^2 \, +\,  {\rm EM}^2 \, ) \,\,-\,\,
\alpha\ib{2}   \, \, \Bigr) \,   $ 

\smallskip \noindent
and the third relation of (3.11) is proven. 

\smallskip \noindent $\bullet \quad$ 
We write now the orthogonality (1.32) between the two edges $\, a =  ({\rm S},\,{\rm
N})  \,$ and  $\, b = ({\rm E},\,{\rm A})  \,$ inside the triangle 
$\, M = ({\rm N},\, {\rm E},\,{\rm A}) .\,$ When we exchange the roles of the  two
vertices $\, {\rm N}\,$ and $\, {\rm E}\,$ in the previous relations, we obtain the
same  result, excepts that $\, \alpha\ib{2} \,$ has to be replaced by $\,\,
\smash{\widetilde{\alpha}}\ib{2} .\,$  So the fourth relation of (3.11) is
established. 

\smallskip \noindent $\bullet \quad$ 
We have from the relation (3.8) : $ \,\,\,  \smash{\widetilde{\alpha}}\ib{2} \,\,=\,\,
{\rm NE}^2 \, \alpha \,\,- \,\, 2 \, {\rm NE}\,\, \alpha\ib{1} \,\,+\,\, \alpha\ib{2}
\,. \,$ Then  

\smallskip \noindent  $   \displaystyle
\alpha\ib{1}  \, \,\,=\,\, {{\bigl( \, {\rm EM}^2 \, \,-\,\, {\rm NM}^2 \, \,+\,\, {\rm
NE}^2 \, \bigr)} \over{2 \, {\rm NE}}} \,  \,\, \alpha \,\,\,=\,\,\, 
 {  \fl {\rm EN} \, \smb \, 
{\bigl( \, \fl{\rm EM}\,+\,  \fl{\rm NM} \, +\,  \fl{\rm EN}\, \bigr)} \over{2 \, 
{\rm NE}}} \,  \,\, \alpha \,$
\smallskip \noindent  $   \displaystyle \qquad
= \,\,  {{  \fl {\rm EN} \, \smb \, \fl{\rm EM} } \over{ \rm NE}}   \,\, \alpha \,$ 

\smallskip \noindent
and the first relation of (3.11) is established. As previously, the exchange of the
vertices $\, {\rm N} \,$ and $ \, {\rm E}\,$ induces the change of $\, \alpha\ib{1} 
\,$ into $\,\, \smash{\widetilde{\alpha}}\ib{1} \,$  that establishes the second
relation of (3.11).

\smallskip \noindent $\bullet \quad$ 
The relations (3.12), (3.13) and (3.14) are
obtained by circular permutation, following the rules that are natural when viewing
the Figure 3 : 

\smallskip \noindent  $   \displaystyle
\alpha \, \longrightarrow \, \beta \, \longrightarrow \, \gamma \, \longrightarrow \,
\delta \,, \quad {\rm E} \, \longrightarrow \, {\rm N} \, \longrightarrow \, {\rm W}
\, \longrightarrow \, {\rm S}  \,, \quad {\rm N} \, \longrightarrow \, {\rm W}
\, \longrightarrow \, {\rm S}  \longrightarrow \, {\rm E} \,$ 

\smallskip \noindent  and $   \displaystyle \quad 
{\rm M} \, \longrightarrow \, {\rm P}  \, \longrightarrow \, {\rm Q}  \,
\longrightarrow \, {\rm R}  \,.\,$ 
$\hfill \square$

\bigskip \noindent
{\bf Lemma 2. \quad A mean value of the dual  Raviart-Thomas basis function}

\noindent
Let  $\, M = ({\rm N},\, {\rm E},\,{\rm A}) \,$ be a triangle of the mesh $\, {\cal
T} \,$ associated to the edge $\, a = ({\rm S},\, {\rm N}) \,$ as in Figure 3 and
$\,\,  \varphi^{\star}\ib{\rm SN} \,\,$ a dual Raviart-Thomas  basis function
satisfying the Hypothesis 2. Then for each constant vector $\, \xi ,\,$ we have : 
\smallskip \noindent  (3.15) $\quad  \displaystyle
\xi \, \smb \, \int_{M} \, \varphi^{\star}\ib{\rm SN} \, \, {\rm d}x \,\,= \,\,
(\xi \, \smb \, n\ib{\rm EN}) \, (\fl{\rm ME} \, \smb \,  n\ib{\rm EN}) \,\,  \alpha
\,.\,$
 
\smallskip 
\vfill \eject 
\noindent
{\bf Proof of Lemma 2}. 

\smallskip \noindent 
We have  : 

\smallskip \noindent   $   \displaystyle  
\xi \, \smb \, \int_{M} \, \varphi^{\star}\ib{\rm SN} \, \, {\rm d}x \,\,= \,\,
\int_{M} \, \nabla \, \bigl(\xi\,\smb\,(x-{\rm M}) \bigr) \,\, \varphi^{\star}\ib{\rm
SN} \, \, {\rm d}x \,  $
\smallskip \noindent   $   \displaystyle    = \,\, 
- \,  \int_{M} \,\, \xi \, \smb \, (x-{\rm M}) \, \, {\rm div}\, 
\varphi^{\star}\ib{\rm SN} \,   \, {\rm d}x \,\,\,+\,\,\,   \int_{\partial M} \,\,
\bigl( \varphi^{\star}\ib{\rm SN} \,\smb\, n\bigr) \,\, \xi \, \smb \, (x-{\rm M}) \,
\, {\rm d}\gamma \,$
\smallskip \noindent   $   \displaystyle    = \,\, 
- \,\bigl( {\rm div}\,  \varphi^{\star}\ib{\rm SN} \bigr) \, \int_{M} \,\, \xi \, \smb
\, (x-{\rm M}) \, \, {\rm d}x \,\,\,+\,\,\,   \int_{\rm NE} \,
\bigl( \varphi^{\star}\ib{\rm SN} \,\smb\, (-n\ib{\rm EN}) \bigr) \,
\xi \, \smb \, (x - {\rm M}) \, \, {\rm d}\gamma \,$
\smallskip \noindent   $   \displaystyle    = \,\, 
0  \,\,\,+\,\,\,   \int_{\rm E N} \, \bigl( \varphi^{\star}\ib{\rm SN} \,\smb\,
(-n\ib{\rm EN}) \bigr) \, \xi \, \smb \, \bigl[ x - {\rm E} + {\rm E} - {\rm M} \bigr]
\, \, {\rm d}\gamma \,$
\smallskip \noindent   $   \displaystyle    = \,\,  
- \int_{\rm E N} \, \bigl( \varphi^{\star}\ib{\rm SN} \,\smb\, n\ib{\rm EN} \bigr) 
\,\xi \, \smb \, \Bigl[ \, s \, {{\fl{\rm EN}}\over{\rm EN}} \,\,+\,\, \fl{\rm ME}
\, \Bigr] \, \, {\rm d}\gamma \,$
\smallskip \noindent   $   \displaystyle    = \,\, 
-\xi \, \smb \,  \Bigl( \,  {{\fl{\rm EN}}\over{\rm EN}} \, \alpha\ib{1} \,\,+\,\, 
\fl{\rm ME}  \, \alpha\,\Bigr) \, \, \,=\,  \,\, 
-\xi \, \smb \, \biggl[ \, \Bigl( \fl{\rm EM}\,\smb\, {{\fl{\rm EN}}\over{\rm EN}}
\Bigr) \,  {{\fl{\rm EN}}\over{\rm EN}} \,\,\,+\,\,\, \fl{\rm ME} \, \biggr] \,\,
\alpha \,\,\,  $
\smallskip \noindent   $   \displaystyle    = \,\, 
-\xi \, \smb \, \Bigl[ \, (\fl{\rm ME} \smb \,   n\ib{\rm EN} ) \,\,   n\ib{\rm EN} \,
\Bigr] \, \,   \alpha \,$ 

\smallskip \noindent 
and the relation (3.15) is established. $ \hfill \square $ 

\bigskip \noindent
{\bf Lemma 3. \quad A simple relation between two triangles}

\noindent
Let  $\, L = ({\rm S},\, {\rm E},\,{\rm N}) \,$ and 
$\, M = ({\rm N},\, {\rm E},\,{\rm A}) \,$ be the two triangles of the mesh $\, {\cal
T} \,$ associated to the edge $\, a = ({\rm S},\, {\rm N}) \,$ as in Figure 3. Let 
$\,\,  \varphi^{\star}\ib{\rm SN} \,\,$ be the  dual Raviart-Thomas  basis function
satisfying the Hypothesis 2. Then  we have : 
\smallskip \noindent  (3.16) $\quad  \displaystyle
{{1}\over{2 \abs{M}}} \,  \int_{M} \,(A-x) \, \smb \,  \varphi^{\star}\ib{\rm SN} \, \,
{\rm d}x \,\,= \,\, {{1}\over{2 \abs{L}}} \, (\fl{\rm SN}\,\smb\,n\ib{\rm EN}) \, 
(\fl{\rm EM} \, \smb\,n\ib{\rm EN}) \,\, \alpha \,.\,$
 
\smallskip \noindent
{\bf Proof of Lemma 3}. 

\smallskip \noindent 
We denote by $\, h \,$ the height of the triangle $\, {\rm NEA} \,$ and by $\, h^* \,$
the height of $\, {\rm ENS}  \,$ chosen such that $\,\,\,  \abs{M} \,=\, {1\over2} \, h
\, {\rm NE} \, \,$ and $\,\,\,  \abs{L} \,=\, {1\over2} \, h^* \, {\rm NE} .\, \,$
We have  from the Lemma 2 :  

\smallskip \noindent   $   \displaystyle    
{{1}\over{2 \abs{M}}} \,  \int_{M} \,(A-x) \, \smb \,  \varphi^{\star}\ib{\rm SN} \,\,
{\rm d}x \,\,= \,\,    
{{1}\over{2 \abs{M}}} \,  \int_{M} \,\bigl( \, N - x + \fl{\rm  NA} \,\bigr)  \, \smb
\,  \varphi^{\star}\ib{\rm SN}  \,{\rm d}x \,$ 
\smallskip \noindent   $   \displaystyle   = \,\, 
{{1}\over{2 \abs{M}}} \,   \fl{\rm NA} \,  \smb \, \int_{M} \, \varphi^{\star}\ib{\rm
SN} \, \, {\rm d}x  \,\,\,= \,\,\, {{1}\over{h \, {\rm NE}}} \,\, (\fl{\rm NA} \, \smb
\, n\ib{\rm EN}) \, (\fl{\rm EM} \, \smb \,  n\ib{\rm EN}) \,\,  \alpha   \,$ 
\smallskip \noindent   $   \displaystyle   = \,\, 
{{1}\over{h \, {\rm NE}}}  \,\, h \,\, (\fl{\rm EM} \, \smb \,  n\ib{\rm EN}) \,\, 
\alpha  \,\,\,= \,\,\, {{1}\over{h^* \, {\rm NE}}}  \,\, h^* \,\, (\fl{\rm EM} \, \smb
\,  n\ib{\rm EN}) \,\,  \alpha \,$
\smallskip \noindent   $   \displaystyle   = \,\, 
{{1}\over{2 \abs{L}}} \, (\fl{\rm SN} \, \smb \,  n\ib{\rm EN}) \,\, (\fl{\rm EM} \,
\smb \,  n\ib{\rm EN}) \,\,  \alpha \,$

\smallskip \noindent 
and the relation (3.16) is proven. $ \hfill \square $ 

\bigskip \noindent
{\bf Proposition 4. \quad Second relation between momenta}

\noindent
The Hypothesis 2  inside the triangles 
$\, L = ({\rm S},\, {\rm E},\,{\rm N}) \,$ and 
$\, M = ({\rm N},\, {\rm E},\,{\rm A}) \,$ implies the following relation : 
\setbox20=\hbox{ $ \displaystyle   
(\rho_{\! L}^2 + {\rm SL}^2) \, \eta \,\,\,+\,\,\, \smash{{1\over6}} \, \fl{\rm SA}
\,\smb\, (\fl{\rm SE} + \fl{\rm SA} +  \fl{\rm SN}) \, \, \alpha \,\,+ \,\,$  }
\setbox21=\hbox{ $ \displaystyle   \qquad  \qquad \qquad \qquad 
\,+ \,\,\, {1\over6} \,  \fl{\rm ND} \,\smb\, (\fl{\rm SD}  \,+\, \fl{\rm SE} \,+\,
\fl{\rm SN}) \,\,\delta \,. \,$  }
\setbox26= \vbox {\halign{# \cr \box20 \cr \box21 \cr  }}
\setbox27= \hbox{ $\vcenter {\box26} $}
\setbox28=\hbox{\noindent (3.17) $\quad \eta\ib{2} \,\,\,= \,\,\, \left\{  \box27
\right. $}
\smallskip  \noindent $ \box28 $
\smallskip  \noindent
We have also : 
\setbox20=\hbox{ $ \displaystyle   
(\rho_{\! L}^2 + {\rm NL}^2) \, \eta \,\,\,+\,\,\, \smash{{1\over6}} \, \fl{\rm SA}
\,\smb\, (\fl{\rm NS} + \fl{\rm NE} +  \fl{\rm NA}) \, \, \alpha \,\,+ \,$  }
\setbox21=\hbox{ $ \displaystyle   \qquad  \qquad \qquad \qquad 
\,+ \,\,\, {1\over6} \,  \fl{\rm ND} \,\smb\, (\fl{\rm NS}  \,+\, \fl{\rm ND} \,+\,
\fl{\rm NE}) \,\,\delta \,\,\,,\,$  }
\setbox26= \vbox {\halign{# \cr \box20 \cr \box21 \cr  }}
\setbox27= \hbox{ $\vcenter {\box26} $}
\setbox28=\hbox{\noindent (3.18) $\quad \widetilde{\eta}\ib{2} \,\,\,= \,\,\, \left\{ 
\box27
\right. $}
\smallskip  \noindent $ \box28 $

\setbox20=\hbox{ $ \displaystyle   
(\rho_{\! K}^2 + {\rm NK}^2) \, \eta \,\,\,-\,\,\, \smash{{1\over6}} \, \fl{\rm SB}
\,\smb\, (\fl{\rm NB} + \fl{\rm NW} +  \fl{\rm NS}) \, \, \beta \,\, \,$  }
\setbox21=\hbox{ $ \displaystyle   \qquad  \qquad \qquad \qquad 
\,- \,\,\, {1\over6} \,  \fl{\rm NC} \,\smb\, (\fl{\rm NW}  \,+\, \fl{\rm NC} \,+\,
\fl{\rm NS}) \,\,\gamma \,\,\,, \,$  }
\setbox26= \vbox {\halign{# \cr \box20 \cr \box21 \cr  }}
\setbox27= \hbox{ $\vcenter {\box26} $}
\setbox28=\hbox{\noindent (3.19) $\quad \widetilde{\eta}\ib{2} \,\,\,= \,\,\, \left\{ 
\box27
\right. $}
\smallskip  \noindent $ \box28 $
\setbox20=\hbox{ $ \displaystyle   
(\rho_{\! K}^2 + {\rm SK}^2) \, \eta \,\,\,-\,\,\, \smash{{1\over6}} \, \fl{\rm SB}
\,\smb\, (\fl{\rm SN} + \fl{\rm SB} +  \fl{\rm SW}) \, \, \beta \,\, \,\,$  }
\setbox21=\hbox{ $ \displaystyle   \qquad  \qquad \qquad \qquad 
\,- \,\,\, {1\over6} \,  \fl{\rm NC} \,\smb\, (\fl{\rm SN}  \,+\, \fl{\rm SW} \,+\,
\fl{\rm SC}) \,\,\gamma \,\,. \,$  }
\setbox26= \vbox {\halign{# \cr \box20 \cr \box21 \cr  }}
\setbox27= \hbox{ $\vcenter {\box26} $}
\setbox28=\hbox{\noindent (3.20) $\quad \eta\ib{2} \,\,\,= \,\,\, \left\{  \box27
\right. $}
\smallskip  \noindent $ \box28 $

\smallskip \noindent
{\bf Proof of Proposition 4}. 

\smallskip \noindent $\bullet \quad$ 
We write the orthogonality (1.32) between the two edges $\, a =  ({\rm S},\,{\rm N}) 
\,$ and  the edge $\, b = ({\rm E},\,{\rm N}) \,$  (see the Figure 3). We have 

\setbox20=\hbox{ $ \displaystyle   
{{1}\over{2 \abs{L}}} \, (x-{\rm S}) \,=\,
{{1}\over{4 \abs{L}}} \,\, \nabla (\abs{ x-{\rm S}}^2) \quad  $
inside    $\, L = ({\rm E},\, {\rm N},\,{\rm S}) \,$  }
\setbox21=\hbox{ $ \displaystyle   
{{1}\over{2 \abs{M}}} \, ({\rm A} - x) \,\,\,\,\,\, \qquad  $ \qquad \qquad \qquad
\quad  inside    $\, M = ({\rm N},\, {\rm E},\,{\rm A}) \,. \,$  }
\setbox26= \vbox {\halign{# \cr \box20 \cr \box21 \cr  }}
\setbox27= \hbox{ $\vcenter {\box26} $}
\setbox28=\hbox{\noindent   $ \varphi\ib{\rm EN}   \,=\, \left\{ 
\box27 \right. $}
\smallskip  \noindent $ \box28 $

\smallskip \noindent 
Then 

\smallskip \noindent  $   \displaystyle
0 \,\,=\,\, \int_{\Omega} \, \varphi^{\star}\ib{\rm SN} \,\smb\,  \varphi\ib{\rm EN} \,
{\rm d}x   \,\,=\,\, 
\int_{L} \, \varphi^{\star}\ib{\rm SN} \,\smb\,  \varphi\ib{\rm
EN} \, {\rm d}x  \,\,+\,\, \int_{M} \, \varphi^{\star}\ib{\rm SN} \,\smb\, 
\varphi\ib{\rm EN} \, {\rm d}x  \,   $ 
\smallskip \noindent  $   \displaystyle
= \,\, {{1}\over{4 \abs{L}}}  \, \int_{L} \, \varphi^{\star}\ib{\rm SN} \,\smb\, 
\nabla (\abs{ x-{\rm S}}^2)  \, {\rm d}x   \,\,\,+\,\,\, {{1}\over{2 \abs{M}}} \, 
 \int_{M} \, \varphi^{\star}\ib{\rm SN} \,\smb\, ( {\rm A} - x) \, {\rm d}x  \,   $ 
\smallskip \noindent  $   \displaystyle
= -\,\, {{1}\over{4 \abs{L}}} \, \bigl( {\rm div}\,  \varphi^{\star}\ib{\rm SN}\bigr)
(L) \,\smb\,  \int_{L} \,\abs{ x-{\rm S}}^2  \, {\rm d}x  \,\,\,+ \,\,  {{1}\over{4
\abs{L}}}  \int_{\partial L} \, (\varphi^{\star}\ib{\rm SN} \,\smb\, n) \,
\abs{ x-{\rm S}}^2    \, {\rm d}\gamma  \, $ 
\smallskip \noindent  $   \displaystyle \qquad  \qquad + \,\, 
{{1}\over{2 \abs{M}}} \,   \int_{M} \, \varphi^{\star}\ib{\rm SN} \,\smb\, ( {\rm A} -
x) \, {\rm d}x  \,$ 
\smallskip \noindent  $   \displaystyle
= -\,\, {{1}\over{4 \abs{L}}} \, \biggl( {{\delta \,+\, \alpha \, - \, \eta
}\over{\abs{L}}} \, \biggr) \,  \int_{L} \,\abs{ x-{\rm S}}^2 \, {\rm d}x  \,\,\,+
\,\,\, $
\smallskip \noindent  $   \displaystyle \qquad  \qquad + \,\, 
{{1}\over{4 \abs{L}}} \, \Bigl( \, \delta\ib{2} \, + \, {\rm SE}^2 \, \alpha \,\,+\,\,
2 \, \fl{\rm SE} \, \smb \, {{ \fl{\rm EN} }\over{\rm EN}} \, \alpha\ib{1} \,\,+\,\,
\alpha\ib{2} \,\,-\,\, \eta\ib{2} \, \Bigr) \,\,+\,\, \, $ 
\smallskip \noindent  $   \displaystyle \qquad  \qquad + \,\, 
{{1}\over{2 \abs{L}}} \,  (\fl{\rm SN} \, \smb \, n\ib{\rm EN}) \, 
(\fl{\rm EM} \, \smb \,  n\ib{\rm EN}) \,\,  \alpha \, \,$

\smallskip \noindent 
and we have, thanks to the relations (3.10), (3.11) and (3.14) :

\setbox20=\hbox{ $ \displaystyle   
-(\delta \,+\, \alpha \, - \, \eta) \,
(\rho\ib{\! L}^2  + {\rm SL}^2) \, \,+\,\, (\rho\ib{\! R}^2  + {\rm SR}^2) \, 
\delta  \,\,+\,\,  {\rm SE}^2 \, \alpha \,\,+\,\, $    }
\setbox21=\hbox{ $ \displaystyle   \quad   + \,\, 
2 \, \Bigl( \fl{\rm SE} \, \smb \, {{ \overrightarrow{\rm EN} }\over{\rm EN}} \Bigr)
\,\Bigl( \fl{\rm EM} \, \smb \, {{ \fl{\rm EN} }\over{\rm EN}} \Bigr) \, \alpha
\,\,+\,\,(\rho\ib{\! M}^2  + {\rm EM}^2) \, \alpha   \,\,+\,\,$ }
\setbox22=\hbox{ $ \displaystyle   \quad   + \,\, 
2 \,  (\fl{\rm SN} \, \smb \, n\ib{\rm EN}) \,  (\overrightarrow{\rm EM} \, \smb \, 
n\ib{\rm EN}) \,\,  \alpha \, \,. \,  $ }
\setbox26= \vbox {\halign{# \cr \box20 \cr \box21 \cr \box22 \cr }}
\setbox27= \hbox{ $\vcenter {\box26} $}
\setbox28=\hbox{\noindent (3.21) $\quad  \eta\ib{2} \,\,=  \left\{  \box27 \right. $}
\smallskip  \noindent $ \box28 $

\smallskip \noindent $\bullet \quad$ 
The $\, \alpha \,$ coefficient   in the right hand side of the relation (3.21) 
 is equal to :  
\smallskip \noindent  $   \displaystyle
-\smash{{{1}\over{36}}}\,({\rm SE}^2 + {\rm EN}^2 + {\rm NS}^2) \,\,-\,\, {\rm SL}^2 
\,\,+\,\,  {\rm SE}^2  \,\,+\,\,  2 \, \bigl( \fl{\rm SE} \, \smb \, \fl{\rm EM} -
(\fl{\rm SE} \, \smb \, n\ib{\rm EN}) \, (\fl{\rm EM} \, \smb \,  n\ib{\rm EN}) \bigr) 
\,\,+   $ 
\smallskip \noindent  $   \displaystyle \qquad  \qquad + \,\, 
{{1}\over{36}}\,({\rm EA}^2 + {\rm AN}^2 + {\rm NE}^2) \,\,+\,\, {\rm EM}^2 
\,\,+\,\,  2 \,  (\fl{\rm SN} \, \smb \, n\ib{\rm EN}) \,  (\fl{\rm EM} \, \smb \, 
n\ib{\rm EN}) \,\,  $
\smallskip \noindent  $   \displaystyle = \,\,\, 
\smash{{{1}\over{36}}}\, \bigl( \, (\fl{\rm EA} + \fl{\rm SE})\,\smb\,(\fl{\rm EA} -
\fl{\rm SE}) \,+\, (\fl{\rm AN} + \fl{\rm NS}) \,\smb\, (\fl{\rm AN} - \fl{\rm NS})
\bigr)  \,\,\,-\,\,\, {\rm SL}^2 + \,\,  (\fl{\rm SE} + \fl{\rm EM})^2 \,$ 
\smallskip \noindent  $   \displaystyle = \qquad 
{{1}\over{36}}\, \fl{ \rm SA} 
\,\smb\,  (\fl{\rm EA} + \fl{\rm ES} + \fl{\rm
NA} + \fl{\rm NS}) \,\,+\,\, (\fl{\rm SM} + \fl{\rm SL}) \,\smb\,  (\fl{\rm SM} -
\fl{\rm SL}) \,$ 
\smallskip \noindent  $   \displaystyle = \,\,\,
{{1}\over{36}}\, \fl{ \rm SA} \,\smb\,   (\fl{\rm EA} + \fl{\rm ES} + \fl{\rm
NA} + \fl{\rm NS}) \,\,+\,\,(2 \, \fl{\rm SL} + \fl{\rm LM}) \,\smb\, \fl{\rm LM} 
\,$ 
\smallskip \noindent  $   \displaystyle = \,\,\,
{{1}\over{36}}\, \fl{ \rm SA} \,\smb\,   (\fl{\rm EA} + \fl{\rm ES} + \fl{\rm
NA} + \fl{\rm NS}) \,\,+\,\,\Bigl( \, {2\over3} \, (\fl{\rm SE} + \fl{\rm SN} ) + 
{1\over3} \,\fl{\rm SA} \, \Bigr) \,\smb\,  \Bigl( \,  {1\over3} \, \fl{\rm SA} \,
\Bigr)  \,$ 
\smallskip \noindent  $   \displaystyle = \,\,\,
{{1}\over{36}}\, \fl{ \rm SA} \,\smb\,  \bigl( \,  \fl{\rm EA} + \fl{\rm ES} + \fl{\rm
NA} + \fl{\rm NS}  + 8 \,( \fl{\rm SE} + \fl{\rm SN} ) + 4 \,\fl{\rm SA} \, \bigr)
\,$ 
\smallskip \noindent  $   \displaystyle = \,\,\,
{{1}\over{36}} \,\fl{ \rm SA} \,\smb\,  \bigl( \, \fl{\rm ES} + \fl{\rm SA} + \fl{\rm
ES} + \fl{\rm NS}  + \fl{\rm SA} +  \fl{\rm NS}  + 8 \,\fl{\rm SE}  + 8 \, \fl{\rm SN}
 + 4 \,\fl{\rm SA} \, \bigr) \,$ 
\smallskip \noindent  $   \displaystyle = \,\,\,
{{1}\over{6}} \, \fl{ \rm SA} \,\smb\,  \bigl( \, \fl{\rm SE}  +  \fl{\rm SA}  + 
\fl{\rm SN}\, \bigr) \,$ 
\smallskip \noindent 
in coherence with the right hand side of the relation (3.17). 

\smallskip \noindent $\bullet \quad$ 
In a similar way, the coefficient of $\, \delta \,$ in the right hand side of the
relation (3.21)   is equal to :  
\smallskip \noindent  $   \displaystyle
-\smash{{{1}\over{36}}}\,({\rm SE}^2 + {\rm EN}^2 + {\rm NS}^2)  
\,\,+\,\, \smash{{{1}\over{36}}}\,( {\rm DE}^2  + {\rm ES}^2  + {\rm SD}^2 ) \,\,+\,\, 
(\fl{\rm SR} + \fl{\rm SL} ) \,\smb\,  (\fl{\rm SR} - \fl{\rm SL} )  \,$ 
\smallskip \noindent  $   \displaystyle = \,\,\,
{{1}\over{36}}\, \bigl( \, (\fl{\rm DE} + \fl{\rm EN }) \,\smb\,    (\fl{\rm DE} -
\fl{\rm EN } ) + (\fl{\rm SD} + \fl{\rm NS}) \,\smb\,   (\fl{\rm SD} - \fl{\rm NS}) \,
\bigr) \,\,\,+ \,$
\smallskip \noindent  $   \displaystyle \qquad  \qquad \qquad  \qquad  + \,\, 
{1\over3} \, ( \fl{\rm SD} + 2 \, \fl{\rm SE} + \fl{\rm SN} ) \,\smb\, 
\Bigl(  {1\over3} \, \fl{\rm ND} \, \Bigr) \,$ 
\smallskip \noindent  $   \displaystyle = \,\,\,
{{1}\over{36}}\, \fl{\rm ND} \,\smb\, ( \, \fl{\rm ED} + \fl{\rm EN} + \fl{\rm SD}  +
\fl{\rm SN} + 4 \, \fl{\rm SD} + 8 \, \fl{\rm SE}  + 4 \, \fl{\rm SN} \, ) \,$ 
\smallskip \noindent  $   \displaystyle = \,\,\,
{{1}\over{36}}\, \fl{\rm ND} \,\smb\, ( \,  \fl{\rm ES} +  \fl{\rm SD} +  \fl{\rm ES} +
\fl{\rm SN} +   \fl{\rm SD}  + \fl{\rm SN} + 4 \, \fl{\rm SD} + 8 \, \fl{\rm SE}  + 4
\, \fl{\rm SN} \, ) \,$ 
\smallskip \noindent  $   \displaystyle = \,\,\,
{{1}\over{6}} \,\fl{\rm ND} \,\smb\, ( \, \fl{\rm SD} +  \fl{\rm SE} + \fl{\rm SN} \, )
\,$ 
\smallskip \noindent 
as proposed in the  right hand side of the relation (3.17). Then the relation (3.17)
is a direct consequence of (3.21).

\smallskip \noindent $\bullet \quad$
The proof of the relation (3.18) is obtained from the previous relation (3.17) with
the following changes :   $ \,\, {\rm E} \longleftrightarrow  {\rm E}  \,,  $
$ \eta \longleftrightarrow \eta \,,  $
$  {\rm A} \longleftrightarrow  {\rm D}  \,,  $ 
$  {\rm N}  \longleftrightarrow  {\rm S}  \,, $ $   M \longleftrightarrow R\,,  $ $
\alpha \longleftrightarrow \delta  \,   $ and $ \,\, \eta\ib{2}
\longleftrightarrow \widetilde{\eta}\ib{2} \,.$ In a similar way, the relations (3.19)
and (3.20) are a straightforward consequence of the relations (3.17) and (3.18) with a
vision of the Figure 3 ``from the top to the bottom'', {\it id est} with the following
changes~:  $  \,\,   {\rm E}  \longleftrightarrow  {\rm W}  \,,  $ $  {\rm N} 
\longleftrightarrow  {\rm S}  \,,  $ $    {\rm D}  \longleftrightarrow  {\rm B}  \,, 
$  $ {\rm A}  \longleftrightarrow  {\rm C}  \,,  $ $   L \longleftrightarrow K \,, 
$ $   M \longleftrightarrow Q \,,  $ $   R \longleftrightarrow P \,,  $ $   \eta
\longleftrightarrow -\eta \,,  $  $   \alpha \longleftrightarrow \gamma \,,  $  
$   \delta \longleftrightarrow \beta \,   $  and 
$ \,\, \eta\ib{2} \longleftrightarrow -\widetilde{\eta}\ib{2} \,.  $
So the proposition is established. $\hfill \square$

\bigskip \noindent
{\bf Proposition 5. \quad Two expressions for the   first order momentum}

\noindent
Under the Hypothesis 2 and the notations proposed at the Figure 3, we have :
\smallskip \noindent  (3.22) $\quad  \displaystyle
\eta\ib{1} \,\,=\,\, \bigl( \, \fl{\rm SL} \, \eta + \fl{\rm LM} \, \alpha + \fl{\rm
LR} \, \delta \, \bigr) \, \smb \, {{\fl{\rm SN}}\over{\rm SN}} \,\,,\,$
\smallskip \noindent  (3.23) $\quad  \displaystyle
\eta\ib{1} \,\,=\,\, \bigl( \, \fl{\rm SK} \, \eta - \fl{\rm KP} \, \beta - \fl{\rm
KQ} \, \gamma \, \bigr) \, \smb \, {{\overrightarrow{\rm SN}}\over{\rm SN}} \,\,.\,$
 
\smallskip \noindent
{\bf Proof of Proposition 5 }. 

\smallskip \noindent $\bullet \quad$ 
We deduce from (3.8), (3.17) and (3.18) : 

\smallskip \noindent  $ \displaystyle
\eta\ib{1} \,\,=\,\, {{1}\over{2 \, {\rm NS}}} \, ( \, \eta\ib{2} -
\widetilde{\eta}\ib{2} + {\rm SN}^2 \, \eta \,) \,$ 
\smallskip \noindent  $ \displaystyle = \,\,\, 
{{1}\over{2 \, {\rm NS}}} \, \Bigl( \, ({\rm SL}^2 - {\rm NL}^2 + {\rm SN}^2 ) \,
\eta \,\,+\,\, \smash{{1\over6}} \, \fl{\rm SA}
\,\smb\, (\fl{\rm SE} + \fl{\rm SA} +  \fl{\rm SN}) \, \, \alpha \,\,+ \,\, $ 
\smallskip \noindent  $ \displaystyle \qquad \qquad + \,\,\, 
{1\over6} \,  \fl{\rm ND} \,\smb\, (\fl{\rm SD}  \,+\, \fl{\rm SE} \,+\,
\fl{\rm SN}) \,\,\delta  \,\
-\,\,\, \smash{{1\over6}} \, \fl{\rm SA}
\,\smb\, (\fl{\rm NS} + \fl{\rm NE} +  \fl{\rm NA}) \, \, \alpha \,$ 
\smallskip \noindent  $ \displaystyle \qquad \qquad - \,\,\, 
 {1\over6} \,  \fl{\rm ND} \,\smb\, (\fl{\rm NS}  \,+\, \fl{\rm ND} \,+\,
\fl{\rm NE}) \,\,\delta \,\Bigr) \,\,\,$  
\smallskip \noindent  $ \displaystyle = \,\,\, 
{{1}\over{2 \, {\rm NS}}} \, \Bigl( \, \bigl[  (\fl{\rm SL} + \fl{\rm NL}) \,\smb\, 
 (\fl{\rm SL} - \fl{\rm NL}) + {\rm SN}^2 \bigr] \,\,  \eta \,$ 
\smallskip \noindent  $ \displaystyle \qquad \qquad + \,\,\, 
{1\over6} \, \fl{\rm SA} \,\smb\, ( \fl{\rm SE} + \fl{\rm SA} +  \fl{\rm SN}  + 
\fl{\rm SN} + \fl{\rm EN}  + \fl{\rm AN}  ) \,\, \alpha \,$ 
\smallskip \noindent  $ \displaystyle \qquad \qquad + \,\,\, 
{1\over6} \, \fl{\rm ND} \,\smb\, ( \fl{\rm SD} + \fl{\rm SE} +  \fl{\rm SN}  + 
\fl{\rm SN} + \fl{\rm DN}  + \fl{\rm EN}  ) \,\, \delta \, \Bigr) \, \,$ 
\smallskip \noindent  $ \displaystyle = \,\,\, 
{{1}\over{2 \, {\rm NS}}} \, \Bigl( \, \fl{\rm SN} \,\smb\, ( \fl{\rm SL} +  \fl{\rm
NL} + \fl{\rm SN} ) \,\, \eta \,\,+\,\, {2\over3} \, \fl{\rm SA} \,\smb\,  \fl{\rm SN}
\,\, \alpha \,\,+\,\, {2\over3} \,  \fl{\rm ND} \,\smb\, \fl{\rm SN}   \,\, \delta \,
\Bigr) \, \,$ 
\smallskip \noindent  $ \displaystyle = \,\,\, 
{{1}\over{{\rm NS}}} \, ( \, \fl{\rm SL} \,\, \eta \,\,+\,\, \fl{\rm LM} \,\,
\alpha \,\,+\,\, \fl{\rm LR} \,\, \delta \,) \,\smb\, \fl{\rm SN}  \,$ 

\smallskip \noindent
and the relation (3.22) is established. 

\smallskip \noindent $\bullet \quad$ 
The proof of the relation (3.23) is analogous. It is a consequence of the relations
(3.8), (3.19) and (3.20) : 

\smallskip \noindent  $ \displaystyle
\eta\ib{1} \,\,=\,\, {{1}\over{2 \, {\rm NS}}} \, ( \, \eta\ib{2} -
\widetilde{\eta}\ib{2} + {\rm SN}^2 \, \eta \,) \,$ 
\smallskip \noindent  $ \displaystyle = \,\,\, 
{{1}\over{2 \, {\rm NS}}} \, \Bigl( \, ({\rm SK}^2 - {\rm NK}^2 + {\rm SN}^2 ) \,
\eta \,\,+\,\, \smash{{1\over6}} \, \fl{\rm SB}
\,\smb\, (\fl{\rm NS} + \fl{\rm BS} +  \fl{\rm WS}) \, \, \beta \,\,+ \,\, $ 
\smallskip \noindent  $ \displaystyle \qquad \qquad + \,\,\, 
{1\over6} \,  \fl{\rm NC} \,\smb\, (\fl{\rm NS}  \,+\, \fl{\rm WS} \,+\,
\fl{\rm CS}) \,\,\gamma  \,\
+\,\,\, \smash{{1\over6}} \, \fl{\rm SB}
\,\smb\, (\fl{\rm NB} + \fl{\rm NW} +  \fl{\rm NS}) \, \, \beta \,\,+ \, $ 
\smallskip \noindent  $ \displaystyle \qquad \qquad + \,\,\, 
 {1\over6} \,  \fl{\rm NC} \,\smb\, (\fl{\rm NW}  \,+\, \fl{\rm NC} \,+\,
\fl{\rm NS}) \,\,\gamma \,\Bigr) \,\,\,$  
\smallskip \noindent  $ \displaystyle = \,\,\, 
{{1}\over{2 \, {\rm NS}}} \, \Bigl( \, \bigl[  (\fl{\rm SK} + \fl{\rm NK}) \,\smb\, 
 (\fl{\rm SK} - \fl{\rm NK}) + {\rm SN}^2 \bigr] \,\,  \eta \,\,+ \,$ 
\smallskip \noindent  $ \displaystyle \qquad \qquad + \,\,\, 
{1\over6} \, \fl{\rm SB} \,\smb\, ( \fl{\rm NS} + \fl{\rm BS} +  \fl{\rm WS}  + 
\fl{\rm NB} + \fl{\rm NW}  + \fl{\rm NS}  ) \,\, \beta \,\,+ \, $ 
\smallskip \noindent  $ \displaystyle \qquad \qquad + \,\,\, 
{1\over6} \, \fl{\rm NC} \,\smb\, ( \fl{\rm NS} + \fl{\rm WS} +  \fl{\rm CS}  + 
\fl{\rm NW} + \fl{\rm NC}  + \fl{\rm NS}  ) \,\, \gamma \, \Bigr) \, \,$ 
\smallskip \noindent  $ \displaystyle = \,\,\, 
{{1}\over{2 \, {\rm NS}}} \, \Bigl( \, \fl{\rm SN} \,\smb\, ( \fl{\rm SK} +  \fl{\rm
NK} + \fl{\rm SN} ) \,\, \eta \,\,+\,\, {2\over3} \, \fl{\rm SB} \,\smb\,  \fl{\rm NS}
\,\, \beta \,\,+\,\, {2\over3} \,  \fl{\rm NC} \,\smb\, \fl{\rm NS}   \,\, \gamma \,
\Bigr) \, \,$ 
\smallskip \noindent  $ \displaystyle = \,\,\, 
{{1}\over{{\rm NS}}} \, ( \, \fl{\rm SK} \,\, \eta \,\,-\,\, \fl{\rm KP} \,\,
\beta \,\,-\,\, \fl{\rm KQ} \,\, \gamma \,) \,\smb\, \fl{\rm SN}  \,. \,$ 

\smallskip \noindent
The relation (3.23) is established and the Proposition 5 is proven.  $\hfill
\square$

\bigskip \noindent
{\bf Lemma 4. \quad Two usefull integrals}

\noindent
Let  $\, K = ({\rm S},\, {\rm N},\,{\rm W}) \,$ and 
$\, L = ({\rm N},\, {\rm S},\,{\rm E}) \,$ be the two triangles of the mesh $\, {\cal
T} \,$ that compose the co-boundary of   the edge $\, a = ({\rm S},\, {\rm N}) \,$ as
in Figure 3.  Let  $\,\,  \varphi^{\star}\ib{\rm SN} \,\,$ be the  a Raviart-Thomas 
basis function satisfying the Hypothesis 2. Then  we have : 
\smallskip \noindent  (3.24) $\quad  \displaystyle
{{1}\over{2 \abs{L}}} \,  \int_{L} \, \varphi^{\star}\ib{\rm SN} \, \smb \, ({\rm
E}-x)  \, \, {\rm d}x \,\,= \,\, {{1}\over{{\rm NS}}} \, ( \, \fl{\rm SL} \,\, \eta
\,\,+\,\,  \fl{\rm LM} \,\, \alpha \,\,+\,\, \fl{\rm LR} \,\, \delta \,) \,\smb\,
n\ib{\rm SN} \,$
\smallskip \noindent  (3.25) $\quad  \displaystyle
{{1}\over{2 \abs{K}}}     \int_{K}  \varphi^{\star}\ib{\rm SN}\, \smb \,  ({\rm W}-x) 
\, \, {\rm d}x \,\,= \,\,{{1}\over{{\rm NS}}} \, ( \fl{\rm NK} \,\, \eta \,\,-\,\, 
\fl{\rm KP} \,  \beta \,\,-\,\,  \fl{\rm KQ} \, \gamma \,) \,\smb\, n\ib{\rm SN} \,.\,$
 
\smallskip \noindent
{\bf Proof of Lemma 4}. 

\smallskip \noindent $\bullet \quad$ 
We establish the relation (3.24) by integrating by parts and using the relations
(3.10), (3.11) and (3.14) :  

\smallskip \noindent $ \displaystyle
{{1}\over{2 \abs{L}}} \,  \int_{L} \, \varphi^{\star}\ib{\rm SN} \, \smb \, ({\rm
E} -x) \, \, {\rm d}x \,\,= \,\, -{{1}\over{4 \abs{L}}} \,  \int_{L} \,
\varphi^{\star}\ib{\rm SN} \, \smb \, \nabla ( \abs{x-{\rm E}}^2 ) \, \, {\rm d}x \,$ 
\smallskip \noindent $ \displaystyle = \,\, 
{{1}\over{4 \abs{L}}} \, ({\rm div}\,  \varphi^{\star}\ib{\rm SN})(L) \,  \int_{L} \, 
\abs{x-{\rm E}}^2 \,  {\rm d}x \,\,-\,\, {{1}\over{4 \abs{L}}} \,  \int_{\partial L}
\, ( \varphi^{\star}\ib{\rm SN} \, \smb \,n) \,  \abs{x-{\rm E}}^2  \, {\rm d}\gamma
\, $ 
\smallskip \noindent $ \displaystyle = \,\, 
{{1}\over{4 \abs{L}}} \,\,\Bigl( \,\, {{ \delta + \alpha - \eta}\over{\abs{L}}} \, 
\int_{L} \,  \abs{x-{\rm E}}^2 \,   {\rm d}x \,\,  $
\smallskip \noindent $ \displaystyle \qquad  \qquad \qquad - \,\Bigl[ \, 
\widetilde{\delta}\ib{2}   \,\,+\,\, \alpha\ib{2}  \,\,-\,\, \bigl( \, \eta\ib{2} 
\,\,+\,\, 2 \, \fl{\rm ES} \,\smb\,  {{\fl{\rm SN}}\over{\rm SN}} \, \eta\ib{1}
\,\,+\,\,  {\rm ES}^2 \,\, \eta \,\, \bigr) \, \Bigr] \, \Bigr) \,$

\smallskip \noindent $ \displaystyle = \,\, 
{{1}\over{4 \abs{L}}} \,\biggl( \, ( \delta + \alpha - \eta ) \, (\rho\ib{\! L}^2
\, +\,  {\rm EL}^2 ) \,\,-\,\, (\rho\ib{\! R}^2 \, +\,  {\rm ER}^2 ) \, \delta
\,\,-\,\, (\rho\ib{\! M}^2 \, +\,  {\rm EM}^2 ) \, \alpha  \,+\,$ 
\smallskip \noindent $ \displaystyle +\,\, 
(\rho_{\! L}^2 + {\rm SL}^2) \, \eta \,\,\,+\,\,\, \smash{{1\over6}} \, \fl{\rm SA}
\,\smb\, (\fl{\rm SE} + \fl{\rm SA} +  \fl{\rm SN}) \, \, \alpha \,+ \,
\smash{1\over6} \,  \fl{\rm ND} \,\smb\, (\fl{\rm SD}  \,+\, \fl{\rm SE} \,+\,
\fl{\rm SN}) \,\delta \, \,$
\smallskip \noindent $ \displaystyle \qquad +\,\, 
2 \, \Bigl( \, \fl{\rm ES} \,\smb\,  {{\overrightarrow{\rm SN}}\over{\rm SN}} \, \Bigr)
\,\,    \bigl( \, \fl{\rm SL} \, \eta + \fl{\rm LM} \, \alpha + \fl{\rm LR} \, \delta
\, \bigr) \,\smb\,  {{\fl{\rm SN}}\over{\rm SN}} \, +\, {\rm ES}^2 \, \eta \, \biggr)
\,$ 

\smallskip \noindent $ \displaystyle = \,\, 
{{1}\over{4 \abs{L}}} \,\Biggl( \, \biggl[\, -{\rm EL}^2 + {\rm SL}^2 + 2 \, 
 \Bigl( \fl{\rm ES} \,\smb\,  {{\overrightarrow{\rm SN}}\over{\rm SN}} \Bigr)  \,
 \Bigl( \fl{\rm SL} \,\smb\,  {{\overrightarrow{\rm SN}}\over{\rm SN}} \Bigr)  
 + {\rm ES}^2 \, \biggl] \,\, \eta \,+\, $ 
\smallskip \noindent $ \displaystyle \quad +\,\, 
\biggl[\,{{1}\over{36}} \, \bigl( (\fl{\rm SE} + \fl{\rm EA}) \,\smb\, (\fl{\rm SE} -
\fl{\rm EA}) + (\fl{\rm SN} + \fl{\rm AN}) \,\smb\, (\fl{\rm SN} - \fl{\rm AN})
\bigr) \,+ \, {\rm EL}^2 -  {\rm EM}^2 \,+ \,$ 
\smallskip \noindent $ \displaystyle \qquad \quad +\,\, 
\smash{{1\over6}} \, \fl{\rm SA} \,\smb\, (\fl{\rm SE} + \fl{\rm SA} +  \fl{\rm SN}) 
\,+\, 2 \, \Bigl(  \fl{\rm ES} \,\smb\,  {{\overrightarrow{\rm SN}}\over{\rm SN}} \,
\Bigr) \, \Bigl( \fl{\rm LM}  \,\smb\,  {{\fl{\rm SN}}\over{\rm SN}} \,
\Bigr) \, \, \biggl] \,\, \alpha \, \,\,+   $
\smallskip \noindent $ \displaystyle \quad +\,\, 
\biggl[\,{{1}\over{36}} \, \bigl( (\fl{\rm NS} + \fl{\rm SD}) \,\smb\, (\fl{\rm NS} -
\fl{\rm SD}) + (\fl{\rm EN} + \fl{\rm DE}) \,\smb\, (\fl{\rm EN} - \fl{\rm DE})
\bigr) \,+ \, {\rm EL}^2 -  {\rm ER}^2 \,+ \,$ 
\smallskip \noindent $ \displaystyle \qquad \quad +\,\, 
\smash{1\over6} \,  \fl{\rm ND} \,\smb\, (\fl{\rm SD} + \fl{\rm SE} + \fl{\rm SN}) \,
\,+\, 2 \, \Bigl(  \fl{\rm ES} \,\smb\,  {{\overrightarrow{\rm SN}}\over{\rm SN}} \,
\Bigr) \, \Bigl( \fl{\rm LR}  \,\smb\,  {{\fl{\rm SN}}\over{\rm SN}} \,
\Bigr) \, \, \biggl] \,\, \delta \,  \Biggr) \,  $

\smallskip \noindent $ \displaystyle = \,\, 
{{1}\over{4 \abs{L}}} \,\Biggl( \, \biggl[\, (\fl{\rm SL} + \fl{\rm EL}) \,\smb\,
(\fl{\rm SL} - \fl{\rm EL}) +  {\rm ES}^2 + 2 \, \fl{\rm ES} \,\smb\, 
\fl{\rm SL} - 2 \,(\fl{\rm ES} \,\smb\, n\ib{\rm SN}) \, (\fl{\rm SL} \,\smb\,
n\ib{\rm SN}) \, \biggl] \,\, \eta \,+\, $ 
\smallskip \noindent $ \displaystyle \qquad \quad +\,\, 
\biggl[\,{{1}\over{36}} \, \fl{\rm SA} \,\smb\, ( \fl{\rm SE}  + \fl{\rm AE}   +
\fl{\rm SN}  + \fl{\rm AN}) \,+\, (\fl{\rm EL} + \fl{\rm EM}) \, (\fl{\rm EL}- 
\fl{\rm EM}) \,+ \,   $
\smallskip \noindent $ \displaystyle \qquad \quad +\,\, 
\smash{{1\over6}} \, \fl{\rm SA} \,\smb\, (\fl{\rm SE} + \fl{\rm SA} +  \fl{\rm SN}) 
\,+\, 2 \,\fl{\rm ES} \,\smb\, \fl{\rm LM} \,-\, 2 \, (\fl{\rm ES} \,\smb\, 
n\ib{\rm SN}) \, (\fl{\rm LM} \,\smb\, n\ib{\rm SN}) \, \biggl] \,\, \alpha \,+\, $ 
\smallskip \noindent $ \displaystyle \qquad \quad +\,\, 
\biggl[\,{{1}\over{36}} \, \fl{\rm ND} \,\smb\, ( \fl{\rm NS}  + \fl{\rm DS}   +
\fl{\rm NE}  + \fl{\rm DE}) \,+\, (\fl{\rm EL} + \fl{\rm ER}) \, (\fl{\rm EL}- 
\fl{\rm ER}) \,+ \,   $
\smallskip \noindent $ \displaystyle \qquad \quad +\,\, 
\smash{1\over6} \,  \fl{\rm ND} \,\smb\, (\fl{\rm SD} + \fl{\rm SE} + \fl{\rm SN}) \,
\,+\, 2 \,\fl{\rm ES} \,\smb\, \fl{\rm LR} \,-\, 2 \, (\fl{\rm ES} \,\smb\, 
n\ib{\rm SN}) \, (\fl{\rm LR} \,\smb\, n\ib{\rm SN}) \, \, \, \biggl] \,\, \delta \, 
\Biggr) \,  $

\smallskip \noindent $ \displaystyle = \,\, 
{{1}\over{4 \abs{L}}} \,\Biggl( \, \biggl[\, \fl{\rm SE} \,\smb\, ( \fl{\rm SL} + 
\fl{\rm EL} + \fl{\rm SE} - 2 \, \fl{\rm SL}) \,+\, 2 \,(\fl{\rm SE} \,\smb\, 
n\ib{\rm SN}) \, (\fl{\rm SL} \,\smb\, n\ib{\rm SN}) \, \biggl] \,\, \eta \,+\, $  
\smallskip \noindent $ \displaystyle \qquad \quad +\,\, 
\biggl[\,{{1}\over{36}} \, \fl{\rm SA} \,\smb\, \bigl(\, \fl{\rm SE}  + \fl{\rm AS}  +
\fl{\rm SE} + \fl{\rm SN} + \fl{\rm AS} + \fl{\rm SN} + 6 \, (\fl{\rm SE} + \fl{\rm
SA} +  \fl{\rm SN}) \, \bigr) \,$
\smallskip \noindent $ \displaystyle \qquad  \qquad \quad +\,\, 
\fl{\rm LM} \,\smb\,( \fl{\rm LE}  + \fl{\rm ME} + 2 \,  \fl{\rm ES} ) \,+\, 
 2 \, (\fl{\rm SE} \,\smb\,  n\ib{\rm SN}) \, (\fl{\rm LM} \,\smb\, n\ib{\rm SN}) \,
\biggl] \,\, \alpha \,+\, $ 
\smallskip \noindent $ \displaystyle \qquad \quad +\,\, 
\biggl[\,{{1}\over{36}} \, \fl{\rm ND} \,\smb\, \bigl(\, \fl{\rm NS}  + \fl{\rm DS}  +
\fl{\rm NS} + \fl{\rm SE} + \fl{\rm DS} + \fl{\rm SE} + 6 \, (\fl{\rm SD} + \fl{\rm
SE} + \fl{\rm SN}) \, \bigr) \,$
\smallskip \noindent $ \displaystyle \qquad  \qquad \quad +\,\, 
\fl{\rm LR} \,\smb\,( \fl{\rm LE}  + \fl{\rm RE} + 2 \,  \fl{\rm ES} ) \,+\, 
 2 \, (\fl{\rm SE} \,\smb\,  n\ib{\rm SN}) \, (\fl{\rm LR} \,\smb\, n\ib{\rm SN}) \, 
 \, \, \biggl] \,\, \delta \,  \Biggr) \,  $
 
\smallskip \noindent $ \displaystyle = \,\, 
{{1}\over{4 \abs{L}}} \,\Bigl( \, \bigl[\,  2 \,(\fl{\rm SE} \,\smb\, 
n\ib{\rm SN}) \, (\fl{\rm SL} \,\smb\, n\ib{\rm SN}) \, \bigl] \,\, \eta \,+\, $  
\smallskip \noindent $ \displaystyle  \quad +\,\, 
\bigl[\,{{1}\over{9}} \, \fl{\rm SA} \,\smb\, ( 2 \, \fl{\rm SE}  + \fl{\rm SA} 
 +  2 \,\fl{\rm SN} ) \,+\, {{1}\over{3}} \, \fl{\rm SA} \,\smb\, ( \fl{\rm LS}  +
\fl{\rm MS}  ) \,+\,   2 \, (\fl{\rm SE} \,\smb\,  n\ib{\rm SN}) \, (\fl{\rm LM}
\,\smb\, n\ib{\rm SN}) \, \bigl] \,\, \alpha \,+\, $  
\smallskip \noindent $ \displaystyle  \quad +\,\, 
\bigl[\,{{1}\over{9}} \, \fl{\rm ND} \,\smb\, ( 2 \, \fl{\rm SE}  + \fl{\rm SN} 
 + \fl{\rm SD} ) \,+\, {{1}\over{3}} \, \fl{\rm ND} \,\smb\, ( \fl{\rm LS}  +
\fl{\rm RS}  ) \,+\,   2 \, (\fl{\rm SE} \,\smb\,  n\ib{\rm SN}) \, (\fl{\rm LR}
\,\smb\, n\ib{\rm SN}) \, \bigl] \,\, \delta \,  \Bigr) \,  $

\smallskip \noindent $ \displaystyle = \,\, 
{{1}\over{2 \abs{L}}} \, (\fl{\rm SE} \,\smb\, n\ib{\rm SN}) \, \Bigl[ \, 
(\fl{\rm SL} \,\smb\, n\ib{\rm SN})   \,\, \eta \,\,+\,\, 
(\fl{\rm LM} \,\smb\, n\ib{\rm SN})   \,\, \alpha \,\,+\,\, 
(\fl{\rm LR} \,\smb\, n\ib{\rm SN})   \,\, \delta \,\Bigr] \,$ 

\smallskip \noindent $ \displaystyle = \,\, 
{{1}\over{\abs {\overrightarrow{\rm NS}}}}  \,  \Bigl[ \, 
(\fl{\rm SL} \,\smb\, n\ib{\rm SN})   \,\, \eta \,\,+\,\, 
(\fl{\rm LM} \,\smb\, n\ib{\rm SN})   \,\, \alpha \,\,+\,\, 
(\fl{\rm LR} \,\smb\, n\ib{\rm SN})   \,\, \delta \,\Bigr] \,$ 

\smallskip \noindent 
that establishes the relation (3.24).

\smallskip \noindent $\bullet \quad$ 
The relation (3.25)  is a consequence of the previous relation (3.24) with the
following modifications~:  $  \,\,   {\rm E}  \longleftrightarrow  {\rm W}  \,,  $ $ 
{\rm N}  \longleftrightarrow  {\rm S}  \,,  $ $    {\rm D}  \longleftrightarrow  {\rm
B}  \,,  $  $ {\rm A}  \longleftrightarrow  {\rm C}  \,,  $ $   L \longleftrightarrow K
\,,  $ $   M \longleftrightarrow Q \,,  $ $   R \longleftrightarrow P \,,  $ $   \eta
\longleftrightarrow -\eta \,,  $  $   \alpha \longleftrightarrow \gamma \,,  $ 
$   \delta \longleftrightarrow \beta \,   $  and 
$ \,\,  n\ib{\rm SN} \longleftrightarrow - n\ib{\rm SN} \,.  $
$\hfill \square$

\bigskip 
\smallskip \noindent {\bf Proof of Theorem  2}. 

\smallskip \noindent $\bullet \quad$ 
We first eliminate the variable $\, \eta\ib{1} \,$ between the relations (3.22) and
(3.23)~; we obtain 

\smallskip \noindent  (3.26) $\quad  \displaystyle
 \bigl( \, \fl{\rm KL} \, \eta \,\,+\,\,  \fl{\rm LM} \, \alpha \,\,+\,\,
 \fl{\rm KP} \, \beta \,\,+\,\, \fl{\rm KQ} \, \gamma \,\,+\,\, \fl{\rm LR} \, \delta 
 \, \bigr) \,\smb\, \fl{\rm SN} \,\,=\,\, 0 \,. \,$ 
 
\smallskip \noindent $\bullet \quad$ 
We write secondly the orthonormality relation (1.32) between the vector function 
$ \,  \varphi\ib{\rm SN} \,$ and its dual $\, \varphi^{\star}\ib{\rm SN} ,\, $ with
the help of (3.24) and (3.25) :  

\smallskip \noindent $ \displaystyle
1 \,\,\, = \,\,\, ( \, \varphi^{\star}\ib{\rm SN}  \,,\, \varphi\ib{\rm SN} \, ) 
\,\,\, = \,\,\,  \int_{K} \,  \varphi^{\star}\ib{\rm SN}  \,\smb\,  \varphi\ib{\rm SN}
\,\, {\rm d}x \,\,\,+\,\,\,  \int_{L} \,  \varphi^{\star}\ib{\rm SN}  \,\smb\, 
\varphi\ib{\rm SN} \,\, {\rm d}x \, \,$
\smallskip \noindent $ \displaystyle = \,\,\, 
{{1}\over{2 \abs{K}}}     \int_{K}  \varphi^{\star}\ib{\rm SN} \,\smb\,  (x-{\rm W}) 
\,\, {\rm d}x \,\,\,+\,\,\,  {{1}\over{2 \abs{L}}} \,  \int_{L} \,
\varphi^{\star}\ib{\rm SN} \, \smb \, ({\rm E}-x)  \, \, {\rm d}x \,\, $ 
\smallskip \noindent $ \displaystyle = \,\,\, 
 {{1}\over{{\rm NS}}} \, \bigl( \, (\fl{\rm SL} +  \fl{\rm KN}) \,\, \eta 
\,\,+\,\,  \fl{\rm LM} \,\, \alpha \,\,+\,\, \fl{\rm KP} \,\,  \beta \,\,+\,\, 
\fl{\rm KQ} \,\, \gamma \,\,+\,\,   \fl{\rm LR} \,\, \delta \, \bigr) \,\smb\,
n\ib{\rm SN} \,$
\smallskip \noindent $ \displaystyle = \,\,\, 
 {{1}\over{{\rm NS}}} \, \bigl( \, (\fl{\rm SL} +  \fl{\rm KN} +  \fl{\rm NS}) \,\,
\eta  \,\,+\,\,  \fl{\rm LM} \,\, \alpha \,\,+\,\, \fl{\rm KP} \,\,  \beta \,\,+\,\, 
\fl{\rm KQ} \,\, \gamma \,\,+\,\,   \fl{\rm LR} \,\, \delta \, \bigr) \,\smb\,
n\ib{\rm SN} \,$
\smallskip \noindent $ \displaystyle = \,\,\, 
 {{1}\over{{\rm NS}}} \, ( \, \fl{\rm KL} \,\, \eta  \,\,+\,\,  \fl{\rm LM} \,\, 
\alpha \,\,+\,\, \fl{\rm KP} \,\,  \beta \,\,+\,\,  \fl{\rm KQ} \,\, \gamma \,\,+\,\,  
\fl{\rm LR} \,\, \delta \, ) \,\smb\, n\ib{\rm SN} \,. \,$

\smallskip \noindent
The relation (3.5) is a direct consequence of the above expression and of the
previous relation (3.26). 

\smallskip \noindent $\bullet \quad$ 
Thirdly, we   eliminate the variable $\, \eta\ib{2} \,$ between the relations (3.17)
and (3.20). The coefficient of the variable $\, \eta \,$ is equal to 

\smallskip \noindent $ \displaystyle 
(\rho_{\! L}^2 + {\rm SL}^2) - (\rho_{\! K}^2 + {\rm SK}^2)  \,\,=\,\, 
{{1}\over{36}} \, \bigl( \, ( {\rm EN}^2 +  {\rm ES}^2 ) - ( {\rm NW}^2 + {\rm SW}^2
)  \, \bigr) \,\,+\,\,  {\rm SL}^2 -  {\rm SK}^2 \,$ 
\smallskip \noindent $ \displaystyle = \, 
{{1}\over{36}} \, \bigl( \, ( \fl{\rm EN} + \fl{\rm NW} ) \smb  ( \fl{\rm EN} +
\fl{\rm WN} ) +  ( \fl{\rm ES} + \fl{\rm SW} )  \smb  ( \fl{\rm ES} + \fl{\rm WS} ) 
 \, \bigr) \,+\,  ( \fl{\rm SL} + \fl{\rm SK} ) \smb  ( \fl{\rm SL} + \fl{\rm
KS} ) \,$ 
\smallskip \noindent $ \displaystyle = \,\,\, 
{{1}\over{36}} \, \fl{\rm WE} \,\smb\, ( \fl{\rm NE} + \fl{\rm NW} +  \fl{\rm SE} +
\fl{\rm SW} ) \,\,+\,\,  {{1}\over{3}} \, \fl{\rm WE} \,\smb\, \biggl( \,
{{1}\over{3}} \, (  \fl{\rm SN} +  \fl{\rm SE } +  \fl{\rm SW} + \fl{\rm SN} ) \,
\biggr) \,$ 
\smallskip \noindent $ \displaystyle = \,\,\, 
{{1}\over{12}} \, \fl{\rm KL} \,\smb\, (\,\fl{\rm NE} + \fl{\rm NW} +  \fl{\rm SE} +
\fl{\rm SW} \,+\, 8 \,  \fl{\rm SN} \,+\, 4 \,  \fl{\rm SE} \,+\, 4 \,  \fl{\rm SW} 
\,) \,$ 
\smallskip \noindent $ \displaystyle = \,\,\, 
{{1}\over{12}} \, \fl{\rm KL} \,\smb\, (\, \fl{\rm NS} + \fl{\rm SE} + \fl{\rm NS} +
\fl{\rm SW} + \fl{\rm SE} + \fl{\rm SW} \,+\, 8 \,  \fl{\rm SN} \,+\, 4 \,  \fl{\rm SE}
\,+\, 4 \,  \fl{\rm SW}  \,) \,$ 
\smallskip \noindent $ \displaystyle = \,\,\, 
{{1}\over{2}} \, \fl{\rm KL} \,\smb\, (\,  \fl{\rm SN} \,+\,    \fl{\rm SE} \,+\,  
\fl{\rm SW}  \,) \,. \,$ 

\smallskip \noindent
We deduce from (3.17), (3.20) and the previous calculus : 

\smallskip \noindent $ \displaystyle 
\fl{\rm KL} \,\smb\, (  \fl{\rm SN} + \fl{\rm SE} + \fl{\rm SW} ) \, \eta \,+\,
\fl{\rm LM} \,\smb\, (  \fl{\rm SE} + \fl{\rm SA} + \fl{\rm SN} ) \,  \alpha \,+\,
\fl{\rm KP} \,\smb\, ( \fl{\rm SN} + \fl{\rm SB} +  \fl{\rm SW})  \, \beta \,+\, $
\smallskip \noindent $ \displaystyle \qquad + \,\, 
\fl{\rm KQ} \,\smb\, ( \fl{\rm SN} + \fl{\rm SW} + \fl{\rm SC})  \,\gamma  \,+\, 
\fl{\rm LR} \,\smb\, ( \fl{\rm SD} + \fl{\rm SE} + \fl{\rm SN}) \, \delta \,\,=\,\, 0
\,$ 

\smallskip \noindent
and taking into consideration the relation  (3.26) : 

\smallskip \noindent $ \displaystyle 
\fl{\rm KL} \,\smb\, (   \fl{\rm SE} + \fl{\rm SW} ) \, \eta \,+\,
\fl{\rm LM} \,\smb\, (  \fl{\rm SE} + \fl{\rm SA}   ) \,  \alpha \,+\,
\fl{\rm KP} \,\smb\, (   \fl{\rm SB} +  \fl{\rm SW})  \, \beta \,+\,  $
\smallskip \noindent $ \displaystyle \qquad + \,\, 
\fl{\rm KQ} \,\smb\, ( \fl{\rm SW} + \fl{\rm SC})  \,\gamma  \,+\, 
\fl{\rm LR} \,\smb\, ( \fl{\rm SD} + \fl{\rm SE}  ) \, \delta \,\,=\,\, 0 \,. \,$ 

\smallskip \noindent
We eliminate the variable $\, \eta \,$ between the previous relation and the relation 
(3.5) after multiplying it by $\,-(\fl{\rm SE} + \fl{\rm SW} ). \,$ We get : 

\smallskip \noindent $ \displaystyle 
\fl{\rm LM} \,\smb\, (  \fl{\rm SE} + \fl{\rm SA} -  \fl{\rm SE} - \fl{\rm SW}  ) \, 
\alpha \,+\, 
\fl{\rm KP} \,\smb\, (   \fl{\rm SB} +  \fl{\rm SW} -  \fl{\rm SE} - \fl{\rm SW} )  \,
\beta \,+\,  $
\smallskip \noindent $ \displaystyle \qquad + \,\, 
\fl{\rm KQ} \,\smb\, ( \fl{\rm SW} + \fl{\rm SC}  -  \fl{\rm SE} - \fl{\rm SW} ) 
\,\gamma  \,+\, 
\fl{\rm LR} \,\smb\, ( \fl{\rm SD} + \fl{\rm SE}  -  \fl{\rm SE} - \fl{\rm SW}  ) \,
\delta \,\,=\,\, $
\smallskip \noindent $ \displaystyle \qquad = \,\, 
- \abs{\rm NS} \,  n\ib{\rm SN}  \,\smb\, ( \fl{\rm SE} + \fl{\rm SW} ) \,$ 

\smallskip \noindent
{\it id est} 

\smallskip \noindent $ \displaystyle 
\fl{\rm LM} \,\smb\, \fl{\rm WA}  \,\, \alpha \,\,+\,\, 
\fl{\rm KP} \,\smb\, \fl{\rm EB}  \,\, \beta  \,\,+\,\,  
\fl{\rm KQ} \,\smb\, \fl{\rm EC}  \,\,\gamma  \,\,+\,\, 
\fl{\rm LR} \,\smb\, \fl{\rm WD}  \,\, \delta \,\,\,\, = \,\,\,  $
\smallskip \noindent $ \displaystyle $ \hfill$  = \,\,\, 
- 3 \, \abs{\rm NS} \,  n\ib{\rm SN}  \,\smb\, ( \fl{\rm SL} + \fl{\rm SK} ) \,$ 

\smallskip \noindent
and the relation (3.6) is a direct consequence of the above relation. The Theorem 2
is established.  $\hfill \square$

\setbox20=\hbox{ $ \displaystyle   
(\rho_{\! L}^2 + {\rm SL}^2) \, \eta \,\,\,+\,\,\, \smash{{1\over6}} \, \fl{\rm SA}
\,\smb\, (\fl{\rm SE} + \fl{\rm SA} +  \fl{\rm SN}) \, \, \alpha \,\,+ \,\,$  }
\setbox21=\hbox{ $ \displaystyle   \qquad  \qquad \qquad \qquad 
\,+ \,\,\, {1\over6} \,  \fl{\rm ND} \,\smb\, (\fl{\rm SD}  \,+\, \fl{\rm SE} \,+\,
\fl{\rm SN}) \,\,\delta \,. \,$  }

\setbox20=\hbox{ $ \displaystyle   
(\rho_{\! K}^2 + {\rm SK}^2) \, \eta \,\,\,-\,\,\, \smash{{1\over6}} \, \fl{\rm SB}
\,\smb\, (\fl{\rm SN} + \fl{\rm SB} +  \fl{\rm SW}) \, \, \beta \,\, \,\,$  }
\setbox21=\hbox{ $ \displaystyle   \qquad  \qquad \qquad \qquad 
\,- \,\,\, {1\over6} \,  \fl{\rm NC} \,\smb\, (\fl{\rm SN}  \,+\, \fl{\rm SW} \,+\,
\fl{\rm SC}) \,\,\gamma \,\,. \,$  }


\bigskip \bigskip
\noindent {\smcaps 4) Perspectives } 

\smallskip \noindent $\bullet \quad$ 
We have proposed to formulate the finite volume method for the Poisson equation in
two space dimensions  with the help of Petrov-Galerkin mixed finite elements. The
unknown is constant in each triangle and the momentum is discretized with the
Raviart-Thomas vectorial finite elements of lower degree. The conservation law is
integrated in each triangle and our   stencil for the discrete gradient operator is
composed by six triangles in the vicinity of each edge  of the mesh. The question of
the determination of such a scheme conducts to a two-parameter family for a possible
choice of a so-called ``dual Raviart-Thomas basis function''  for the finite volume
scheme.  We have also developed a sufficient hypothesis to prove  the stability and the
optimal convergence of the associated finite volume scheme.  The next step of this
research is to construct explicitly an interpolation vector valued function  in the
particular case where $\, \Omega = \R^2 \,$ in order to  determine  free coefficients
and to establish the stability property.  

\noindent $\bullet \quad$ 
We thank Jean-Pierre Croisille for his kind invitation to first present  [Du99]  the
 results contained in this 	article and for  regular helpfull discussions  that 
 convinced us of the complexity of  mathematical links   between the present
formulation of the  finite volume method   and his analysis [Cr2k] of the box scheme.


\bigskip \bigskip  \noindent {\smcaps 5) References}

\smallskip \hangindent=9mm \hangafter=1 \noindent  [Ba71]  $\,\,$  
 I. Babu\v{s}ka,   Error-bounds for finite
element method, {\it Numerische Mathematik}, vol.~16, p. 322-333, 1971. 

\smallskip \hangindent=9mm \hangafter=1 \noindent  [CVV99]   $\,\,$ 
 Y. Coudi\`ere, J.P. Vila, P. Villedieu,
Convergence rate of a finite volume scheme for a two-dimensional convection-diffusion
problem, {\it Mathematical Modelling and Numerical Analysis}, vol. 33, p. 493-516,
1999.  

\smallskip \hangindent=9mm \hangafter=1 \noindent  [BDM85]  $\,\,$ 
 F. Brezzi, J. Douglas, L.D. Marini, 
 Two families of mixed finite elements for second order elliptic problems, 
{\it Numerische Mathematik}, vol. 47, p. 217-235, 1985. 

\smallskip \hangindent=9mm \hangafter=1 \noindent  [BMO96]  $\,\,$  
J. Baranger, J.F. Ma\^{\i}tre, F. Oudin,
Connection between finite volumes and mixed finite element methods, {\it Mathematical
Modelling and Numerical Analysis}, vol. 30, p. 445-465, 1996.

\smallskip \hangindent=9mm \hangafter=1 \noindent  [CR72]  $\,\,$  
 P.G. Ciarlet, P.A. Raviart,   General
Lagrange and  Hermite interpolation in $R^n$ with applications to finite element
methods, {\it  Archive for Rational Mechanics and Analysis}, vol. 46, p.177-199,
1972.

\smallskip \hangindent=9mm \hangafter=1 \noindent  [Cr2k] $\,\,$  
  J.C. Croisille, Finite Volume Box Schemes
and Mixed Methods, {\it Mathematical Modelling and Numerical Analysis}, vol. 34, n$^o$
5, p. 1087-1106, 2000. 

\smallskip \hangindent=9mm \hangafter=1 \noindent  [Du89]  $\,\,$  
 F. Dubois, Calcul des flux visqueux dans un
code de r\'esolution des \'equations de Navier Stokes par une m\'ethode de volumes
finis non structur\'es,  {\it Aerospatiale Les Mureaux}, internal  report SDTMI
104/89, july 1989.

\smallskip \hangindent=9mm \hangafter=1 \noindent  [Du92] $\,\,$  
F. Dubois, Interpolation de Lagrange et
volumes finis, {\it Aerospatiale Les Mureaux}, internal  report  STS 104109, february
1992, see also {\it Lemmes finis pour la dynamique des gaz}, to appear.

\smallskip \hangindent=9mm \hangafter=1 \noindent  [Du99] $\,\,$ 
 F. Dubois, Volumes finis par \'el\'ements
finis mixtes de Petrov-Galerkin~: vers le bidimensionnel, {\it Seminar, Metz
University}, october 21, 1999.

\smallskip \hangindent=9mm \hangafter=1 \noindent  [Du2k] $\,\,$  
 F. Dubois, Finite Volumes and Mixed
Petrov-Galerkin Finite Elements: The Unidimensional Problem,   {\it Numer.
Meth. in Part. Diff. Eq.}, vol. 16, n$^{\rm o}\,$3, p. 335-360,  2000.

\smallskip \hangindent=9mm \hangafter=1 \noindent   [FGH91] $\,\,$  
 I. Faille, T. Gallou\"et, R. Herbin, Les
math\'ematiciens d\'ecouvrent les volumes finis, {\it Matapli}, vol. 23, p. 37-48,
octobre 1991.

\smallskip \hangindent=9mm \hangafter=1 \noindent    [Go71]  $\,\,$  
C. Godbillon, {\it Elements de topologie
alg\'ebrique}, Hermann, Paris, 1971. 

\smallskip \hangindent=9mm \hangafter=1 \noindent  [He95]   $\,\,$  
R. Herbin, An error estimate for a finite
volume scheme for a diffusion-convection problem in a triangular mesh, {\it Numer.
Meth. in Part. Diff. Eq.}, vol. 11, p. 165-173, 1995.

\smallskip \hangindent=9mm \hangafter=1 \noindent   [Ke81]  $\,\,$  
D.S. Kershaw,  Differencing of the Diffusion
Equation in Lagrangian Hydrodynamic Codes,  {\it J. of Computational Physics}, vol. 39,
p. 375-395, 1981.

\smallskip \hangindent=9mm \hangafter=1 \noindent  [No64]  $\,\,$  
W.F. Noh, ``CEL: a  time dependent two space
dimensional, coupled Euler-Lagrange code'', {\it Methods in computational physics}, 
vol. 3, p. 117-179, Academic Press, New York, 1964.

\smallskip \hangindent=9mm \hangafter=1 \noindent [Pa80]  $\,\,$  
 S.V. Patankar,  {\it Numerical Heat Transfer
and Fluid Flow}, Hemisphere, 1980. 

\smallskip \hangindent=9mm \hangafter=1 \noindent   [RT77]  $\,\,$  
P.A. Raviart, J.M. Thomas, ``A mixed finite
element method for 2nd order elliptic problems'', in {\it Lecture Notes in
Mathematics}, vol. 606 (Dold-Eckmann Eds), Springer-Verlag, Berlin, p. 292-315, 1977.

\smallskip \hangindent=9mm \hangafter=1 \noindent   [RT91]  $\,\,$ 
J.E. Roberts J.M. Thomas, Mixed and Hybrid
Methods, {\it Handbook of Numerical Analysis} (Ciarlet-Lions Eds),  vol. II, Finite
Element Methods (Part I), p. 523-639, Elsevier Science Publishers, Amsterdan,
1991.

\smallskip \hangindent=9mm \hangafter=1 \noindent  [TT99]  $\,\,$ 
 J.M. Thomas, D. Trujillo, Mixed finite
volume methods, {\it Int. J. Numer. Meth. Eng.}, vol. 46, p. 1351-1366, 1999. 

\smallskip \hangindent=9mm \hangafter=1 \noindent [YMAC97]  $\,\,$ 
A. Youn\`es, R. Mose, P. Ackerer, G.
Chavent, Une r\'esolution par les \'el\'ements finis mixtes \`a une inconnue par
maille,  {\it Comptes Rendus Acad. Sci. Paris}, S\'erie I, vol. 328, n$^o$7, p.
623-626, 1997.


\bye